\documentclass[a4paper,10pt,psamsfonts]{amsart}
\usepackage[utf8x]{inputenc}
\usepackage[T1]{fontenc}
\usepackage[ngerman,russian,american,english]{babel}

\usepackage[pdftex]{graphicx}

\usepackage{textcomp}

\usepackage{comment}

\DeclareFontFamily{U}{matha}{\hyphenchar\font45}
\DeclareFontShape{U}{matha}{m}{n}{
      <5> <6> <7> <8> <9> <10> gen * matha
      <10.95> matha10 <12> <14.4> <17.28> <20.74> <24.88> matha12
      }{}
\DeclareSymbolFont{matha}{U}{matha}{m}{n}
\DeclareFontFamily{U}{mathx}{\hyphenchar\font45}
\DeclareFontShape{U}{mathx}{m}{n}{
      <5> <6> <7> <8> <9> <10>
      <10.95> <12> <14.4> <17.28> <20.74> <24.88>
      mathx10
      }{}
\DeclareSymbolFont{mathx}{U}{mathx}{m}{n}
\DeclareMathSymbol{\obot}         {2}{matha}{"6B}
\DeclareMathSymbol{\bigobot}       {1}{mathx}{"CB}

\usepackage{mathrsfs}

\usepackage{amsfonts}
\usepackage{amsmath}
\usepackage{amssymb}
\usepackage{amsthm}
\usepackage{stmaryrd}

\usepackage{algorithm}
\usepackage{algorithmic}

\usepackage{enumerate}

\usepackage{amscd}
\usepackage{tikz}
\usepackage{tikz-cd}
\usepackage{sidecap}
\usetikzlibrary{matrix,arrows}

\usepackage{listings}  

\usepackage{float}
\floatstyle{boxed}
\restylefloat{table}
\restylefloat{figure}

\usepackage{todonotes}

\usepackage{framed}

\newtheoremstyle{break}{9pt}{9pt}{}{}{\bfseries}{.}{\newline}{}
\theoremstyle{break}

\newtheorem{theorem}{Theorem}[section]
\newtheorem{lemma}[theorem]{Lemma}

\newtheorem{remark}[theorem]{Remark}
\newtheorem{example}[theorem]{Example}

\newtheorem{definition}[theorem]{Definition}

\newtheorem{corollary}[theorem]{Corollary}

\newcommand{\comments}[1]{}

\newcommand\explanationend{\xqed{$\bigstar\surd$}}

\newcommand{\vol}{\operatorname{vol}}

\newcommand{\inv}{{-1}}

\newcommand{\dif}{{\mathrm d}}

\newcommand{\grad}{\operatorname{grad}}

\newcommand{\curl}{\operatorname{curl}}

\newcommand{\rot}{\curl}

\newcommand{\diver}{\operatorname{div}}

\newcommand{\st}{ \mid }

\newcommand{\Kappa}{{\Gamma}}

\newcommand{\trace}{\operatorname{tr}}
\newcommand{\normaltrace}{\operatorname{nm}}

\newcommand{\ext}{\operatorname{ext}}

\newcommand{\Id}{\operatorname{Id}}

\newcommand{\dom}{\operatorname{dom}}

\newcommand{\rng}{\operatorname{ran}}

\newcommand{\cartan}{{\mathsf d}}

\newcommand{\cocartan}{{\mathsf{\delta}}}

\newcommand{\Difftot}{{\cartan}}
\newcommand{\Diffhor}{{\sfD}}
\newcommand{\Diffver}{{\sfT}}
\newcommand{\Antihor}{{\sfP}}
\newcommand{\Antiver}{{\sfE}}

\newcommand{\cellless}{{ {} \vartriangleleft {} }}

\newcommand{\cellleq}{{ {} \trianglelefteq {} }}

\newcommand{\ncellleq}{{ {} \ntrianglelefteq {} }}

\newcommand\xqed[1]{%
  \leavevmode\unskip\penalty9999 \hbox{}\nobreak\hfill
  \quad\hbox{#1}}
\newcommand\demo{\xqed{$\triangle$}}

\newcommand\xqedhere[2]{
  \rlap{\hbox to#1{\hfil\llap{\ensuremath{#2}}}}}

\newcommand{\bbN}{{\mathbb N}}

\newcommand{\bbR}{{\mathbb R}}

\newcommand{\calC}{{\mathcal C}}

\newcommand{\calH}{{\mathcal H}}

\newcommand{\calM}{{\mathcal M}}
\newcommand{\calN}{{\mathcal N}}

\newcommand{\calP}{{\mathcal P}}

\newcommand{\calR}{{\mathcal R}}

\newcommand{\calT}{{\mathcal T}}
\newcommand{\calU}{{\mathcal U}}
\newcommand{\calV}{{\mathcal V}}

\newcommand{\frakC}{{\mathfrak C}}

\newcommand{\frakH}{{\mathfrak H}}

\newcommand{\sfD}{{\mathsf D}}
\newcommand{\sfE}{{\mathsf E}}

\newcommand{\sfP}{{\mathsf P}}

\newcommand{\sfT}{{\mathsf T}}

\lstdefinelanguage{pseudocode}
{
  morekeywords = {
    if, then, else,
    select, case,
    for, each, in, from, to, next, break,
    goto,
    let,
    do, while, until,
    return
    function, input, output,
    asdf
  },
  sensitive=false,
  morecomment=[l]{--},
  morecomment=[l]{//},
  morecomment=[s]{/*}{*/},
  morestring=[b]",
}

\title[Homology theory of DDDF]
{Complexes of Discrete Distributional Differential Forms 
and their Homology Theory}

\author{Martin Werner Licht}

\address{Department of Mathematics, University of Oslo, PO box 1053 Blindern, NO-0316 Oslo, Norway}

\email{martinwl@math.uio.no}

\thanks{This research was supported by the European Research Council through the FP7-IDEAS-ERC Starting Grant scheme, project 278011 STUCCOFIELDS.}



\allowdisplaybreaks

\begin{document}

\def\justbeingincluded{justbeingincluded}



\begin{abstract}
 Complexes of discrete distributional differential forms are introduced into finite element exterior calculus.
 Thus we generalize a notion of Braess and Sch\"oberl, originally studied for a posteriori error estimation.
 We construct isomorphisms between the simplicial homology groups of the triangulation,
 the discrete harmonic forms of the finite element complex,
 and the harmonic forms of the distributional finite element complexes.
 As an application, we prove that the complexes of finite element exterior calculus
 have cohomology groups isomorphic to the de Rham cohomology, including the case of partial boundary conditions.
 Poincar\'e-Friedrichs-type inequalities will be studied in a subsequent contribution.
\end{abstract}

\maketitle


\ifx\justbeingincluded\undefined
\input{../global/header_common.tex}


\begin{document}


\fi

\section{Introduction}


Finite element exterior calculus (FEEC, \cite{AFW1,AFW2}) recasts finite element theory 
in the calculus of differential forms and has emerged as a unifying framework for vector-valued finite elements.
The classical residual error estimator has been studied recently 
in the setting of finite element exterior calculus by Demlow and Hirani \cite{demlow2014posteriori}.
On the contrary, implicit error estimators have generally remained focused to spaces of scalar functions
in literature (see \cite{ainsworth2011posteriori,repin2008posteriori,verfurth2013posteriori}),
despite their promising performance in numerical experiments \cite{carstensen2010estimator}.
A notable exception is the equilibrated residual error estimator 
that Braess and Schöberl \cite{BrSchoMax} have studied for first-order N\'ed\'elec elements in two and three dimensions.

They have introduced distributional finite element complexes
--- the major contribution of the present work is to integrate that notion into finite element exterior calculus.
We formulate complexes of discrete distributional differential forms and determine their homology spaces.
A subsequent work will analyze Poincar\'e-Friedrichs inequalities of these complexes.
The present work also serves as a technical preparation for research on a posteriori error estimation within FEEC,
but we derive results of independent interest.
For example, we close a gap in literature and derive compatibility on homology
for the standard finite element complex in the case of partial boundary conditions,
which is relevant for the Hodge Laplace equation with mixed boundary conditions.
\\



We outline the essential ideas with an example --- here we employ the formalism of vector calculus, close to \cite{BrSchoMax},
but the remainder of the contribution employs the calculus of differential forms.
Let $\Omega \subset \bbR^{3}$ be a Lipschitz domain with a triangulation $\calT$.
We denote by $\calT^3$, $\calT^2$, $\calT^1$ and $\calT^0$ the sets
of tetrahedrons, triangles, edges and vertices of $\calT$, respectively.

With respect to this triangulation, we consider a first-order finite element complex:
\begin{align}
 \label{eq:intro_fem_sequence}
 \begin{CD}
 \calP^1(\calT)_0 @>{\grad}>> \mathcal{N}d(\calT)_0 @>{\rot}>> \mathcal{R}T(\calT)_0 @>{\diver}>> \calP_{-1}^0(\calT).
 \end{CD}
\end{align}
Here,
$\calP^1(\calT)_0$ is the space of continuous piecewise affine functions satisfying homogeneous boundary conditions,
$\mathcal{N}d(\calT)_0$ is the curl-conforming first-order N\'ed\'elec space satisfying homogeneous tangential boundary conditions,
$\mathcal{R}T(\calT)_0$ is the divergence-conforming first-order Raviart-Thomas space satisfying homogeneous normal boundary conditions,
and $\calP_{-1}^0(\calT)$ is the space of piecewise constant functions.

Following the notation of \cite{BrSchoMax},
we let $\calP_{-1}^1(\calT)$ be the space of piecewise affine functions;
loosely speaking, this is $\calP^1(\calT)$ without the requirement of boundary conditions and continuity along faces.
Although the classical gradient is not defined on $\calP_{-1}^1(\calT)$,
we can view $\calP_{-1}^1(\calT)$ as a space of functionals on smooth functions
and then apply the gradient in the sense of distributions.
The distributional gradient of $u \in \calP_{-1}^1(\calT)$ is defined via
\begin{align}
 \label{eq:distributional_gradient_definition}
 \left\langle\; \grad u, \vec v \;\right\rangle
 :=
 - \sum_{T \in \calT^3} \int_T u \diver \vec v \;\dif x,
 \quad
 \vec v \in C^\infty(\overline\Omega,\bbR^3).
\end{align}
Then piecewise application of the divergence theorem shows:
\begin{align}
 \label{eq:divergence_theorem_application}
 \left\langle\; \grad u, \vec v \;\right\rangle
 =
 \sum_{T \in \calT^3} \int_T \langle\; \grad u, \vec v \;\rangle \;\dif x
 -
 \sum_{T \in \calT^3} \sum_{\substack{F \in \calT^2\\F \subset T}}
 \int_{F} u \langle \vec v, \vec n_{T,F} \rangle \;\dif s,
\end{align}
where $\vec n_{T,F}$ is the outward normal along the face $F$ of $T$.
The distributional gradient maps $\calP_{-1}^1(\calT)$ into a space of functionals on $C^\infty(\overline\Omega,\bbR^3)$,
which we denote by $\mathcal{N}d_{-2}(\calT)$, and which is spanned by two types of functionals: 
on the one hand, by integration against piecewise first-order N\'ed\'elec elements, 
which are not necessarily tangentially continuous along faces,
and on the other hand, by functionals that act as integration against affine functions over faces.
This corresponds to the two terms on right-hand side of \eqref{eq:divergence_theorem_application}.

Next, the $\curl$-operator, in the sense of distributions, maps $\mathcal{N}d_{-2}(\calT)$ 
into another space $\mathcal{R}T_{-3}(\calT)$ of functionals on $C^\infty(\overline\Omega,\bbR^3)$.
The space $\mathcal{R}T_{-3}(\calT)$ is spanned by integrals against piecewise first-order Raviart-Thomas elements,
which are not necessarily normally continuous along faces,
by integral functionals along faces, and integral functionals along edges.
Eventually, the distributional divergence maps $\mathcal{R}T_{-3}(\calT)$ into $\calP^0_{-4}(\calT)$,
which is a space of functionals over $C^\infty(\overline\Omega)$ again,
spanned by integral evaluations on tetrahedrons, faces and edges, and by point evaluations.
We have thus found a distributional finite element complex \cite[Equation (3.18)]{BrSchoMax}:
\begin{align}
 \label{eq:intro_distributional_sequence}
 \begin{CD}
 \calP^1_{-1}(\calT) @>{\grad}>> \mathcal{N}d_{-2}(\calT) @>{\rot}>> \mathcal{R}T_{-3}(\calT) @>{\diver}>> \calP^0_{-4}(\calT).
 \end{CD}
\end{align}
See also \cite[Equations (3.3), (3.5), (3.7), (3.16-3.18)]{BrSchoMax} for similar differential complexes.
The publication \cite{BrSchoMax} of Braess and Schöberl considers distributional finite element complexes
on local patches of two- and three-dimensional triangulations,
based on finite element spaces of first order.
A unified and more general treatment of these distributional finite element complexes
is possible with the calculus of differential forms.
Thus we extend the basic idea to discrete distributional de Rham complexes of any dimension, 
over domains of arbitrary topology, and with general partial boundary conditions.
The theory includes the spaces of piecewise polynomial differential forms of finite element exterior calculus.


Distributional differential forms and similar ideas appear in different areas of mathematics.
De Rham \cite{rham1984differentiable} introduced the term ``currents'' for continuous linear functionals
on a class of locally convex spaces of smooth differential forms.
Geometric integration theory \cite{krantz2008geometric} knows simplicial chain complexes as a specific example of currents,
which is also rediscovered in this work.
Christiansen \cite{christiansen2011linearization} has considered distributional finite element complexes in Regge calculus.
\\



We can view the right-hand side of \eqref{eq:divergence_theorem_application} as the sum of two operators:
a piecewise differential operator on the one hand, and a ``jump term'' that is the sum of signed traces on the other hand.
Similar decompositions hold for the distributional curl and the distributional divergence,
and a unified treatment is accessible with the calculus of differential forms.

It is an essential observation of this contribution that both of these operators are constituent for differential complexes.
For instance, the standard finite element complex \eqref{eq:intro_fem_sequence} is a subcomplex of the distributional finite element complex \eqref{eq:intro_distributional_sequence},
and it composed of spaces on which the last term of \eqref{eq:divergence_theorem_application} vanishes.
But when identifying simplices with their indicator functions,
we furthermore observe that, up to a sign convention, the simplicial chain complex of the triangulation
\begin{align}
 \label{seq:intro_chain_sequence}
 \begin{CD}
  \calC_3(\calT) @>{-\partial^{3}}>> \calC_2(\calT) @>{\partial^{2}}>> \calC_1(\calT) @>{-\partial^{1}}>> \calC_0(\calT)
 \end{CD}
\end{align}
is a subcomplex of \eqref{eq:intro_distributional_sequence} as well.
It is composed of spaces on which the piecewise differential vanishes.
Indeed, a close look reveals that the ``jump term`` in the right-hand side of \eqref{eq:divergence_theorem_application}
resembles the simplicial boundary operator.
Given these differential complexes, the incentive question that we address in this contribution is:
can we relate their homology spaces?

To answer this question, it we employ the notion of \emph{double complex} that is well-known in homological algebra.
It is instructive to consider the following diagram:
\begin{align}
 \label{eq:intro_doublecomplex}
 \begin{CD}
  \calP^1_{-1}(\calT^3) @>{ \grad_{\calT}}>> \calN d_{-1}(\calT^3) @>{ \curl_{\calT}}>> \calR T_{-1}(\calT^3) @>{ \diver_{\calT}}>> \calP^0_{-1}(\calT^3)
  \\
  @V{-\partial^{3}}VV  @V{-\partial^{3}}VV  @V{-\partial^{3}}VV  @. 
  \\
  \calP^1_{-1}(\calT^2) @>{-\grad_{\calT}}>> \calN d_{-1}(\calT^2) @>{-\curl_{\calT}}>> \calR T_{-1}(\calT^2)
  \\
  @V{\partial^{2}}VV  @V{\partial^{2}}VV  @.   @. 
  \\
  \calP^1_{-1}(\calT^1) @>\grad_{\calT}>> \calN d_{-1}(\calT^1) @. {}
  \\
  @V{-\partial^{1}}VV  @.  @.   @. 
  \\
  \calP^1_{-1}(\calT^0) @. {} @. {}  
 \end{CD}
\end{align}
The spaces in this diagram are finite element spaces on lower-dimensional skeletons of the triangulation
without continuity conditions imposed.
For instance, 
$\calN d_{-1}(\calT^2)$ is a space of functionals on $C^{\infty}(\overline\Omega,\bbR^{3})$
which act by integration against first-order N\'ed\'elec elements over two-dimensional simplices,
and
$\calP^1_{-1}(\calT^0)$ is the span of point evaluations at vertices of $\calT$ acting on $C^{\infty}(\overline\Omega)$;
similarly for the other spaces in \eqref{eq:intro_doublecomplex}.
The horizontal mappings are piecewise differential operators, and thus the rows of the diagram are differential complexes by themselves.
The rows are furthermore exact sequences for many choices of finite spaces,
including the exact sequences of finite element exterior calculus.
The vertical mappings correspond to the boundary terms in partial integration formulas like \eqref{eq:divergence_theorem_application},
and in the above diagram we suggestively denote them by the same symbol as the simplicial boundary operator.
It is an original observation of this work that these ''jump-terms`` are operators in their own right
and constituent for differential complexes in the columns.
The columns are even exact sequences;
this uses the geometric decomposition of the finite element spaces
and a combinatorial condition on the triangulation.
The diagram \eqref{eq:intro_doublecomplex} is a \emph{double complex} in the sense of homological algebra \cite{gelfand1999homological},
and we note that the complex \eqref{eq:intro_distributional_sequence} corresponds to the sequence of diagonals,
also called \emph{total complex}, of the double complex.
Furthermore, our earlier observations transfer: the simplicial chain complex is included in the left-most column,
whereas the standard finite element complex is included in the top-most row.

The question regarding the homology spaces can be answered with the adaption of methods that evolved in the treatment of double complexes.
We construct isomorphisms between the homology groups of the triangulation,
the discrete harmonic forms of the standard finite element complex,
and discrete distributional harmonic forms of distributional finite element complexes such as \eqref{eq:intro_distributional_sequence}.
This is an alternative access towards homology theory in a finite element setting,
besides de Rham mappings \cite{StructPresDisc} and smoothed projections \cite{christiansen2008smoothed},
and, to the author's best knowledge, this is first derivation in literature of the homology theory of finite element de Rham complexes with partial boundary conditions.


Double complexes are used in differential topology,
with the \v Cech de Rham complex being the most prominent example \cite{weil1952theoremes}. 
Falk and Winther \cite{falk2012local} have recently introduced a finite element \v Cech de Rham complex to finite element theory,
albeit not for questions of homological nature.
\\


The remainder of this contribution is organized as follows.
In Section \ref{sect:preparation} we recall relevant notions on simplicial complexes, differential forms, and Hilbert complexes.
Furthermore, the $L^2$ de Rham complex is introduced as an example application of this contribution's theory.
In Section \ref{sect:element_systems} we introduce discrete distributional differential forms and relevant differential operators.
In Section \ref{sect:homology_theory_theoretical} we study finite element double complexes such as \eqref{eq:intro_doublecomplex}.
We prove the exactness of the rows and columns under reasonable assumptions.
This shows the existence of an isomorphism between the simplicial homology groups and the discrete harmonic forms.
In Section \ref{sect:distributional_homology} we study complexes of discrete distributional differential forms such as \eqref{eq:intro_distributional_sequence},
and derive isomorphisms between the discrete distributional harmonic forms.
This provides also a constructive proof of the result in the preceding section.

\ifx\justbeingincluded\undefined
\end{document}
\fi


\ifx\justbeingincluded\undefined
\input{../global/header_common.tex}

\begin{document}

\fi

\section{Preliminaries}
\label{sect:preparation}

In this section we gather technical prerequisites and notational conventions from topology,
analysis on manifolds, and functional analysis.
This work focuses on finite element spaces on simplices, differential and trace operators between them,
and the homology of finite-dimensional Hilbert complexes. 
Therefore we review simplicial complexes (Subsection \ref{subsect:preparation_simplicial}), 
differential forms on manifolds (Subsection \ref{subsect:preparation_differentialforms}), 
and basic aspects of Hilbert complexes (Subsection \ref{subsect:preparation_hilbertcomplexes}).
We also give an outline of the $L^2$ de Rham complex of polyhedrally bounded domains with partial boundary conditions
(Subsection \ref{subsect:preparation_derhamcomplex})
in order to show the connections to the analysis of partial differential equations,
and in order to introduce an example application that is repeatedly addressed in the sequel.

\subsection{Simplicial complexes}
\label{subsect:preparation_simplicial}

We review basic notions of simplicial topology and the homology of simplicial chain complexes.
We refer to the textbooks \cite{LeeTopological}, \cite{mac1975homology} and \cite{spanier1995algebraic}
for further background on these topics.
\\

Let $n \in \bbN_{0}$ be fixed.
A (closed) \emph{$m$-simplex} $C$ is the convex closure of a set of $m+1$ affinely independent points in $\bbR^n$,
which we call the \emph{vertices} of $C$, and we then also write $\dim C = m$.
A simplex $F$ is called a \emph{subsimplex} of an $m$-simplex $C$
if the set of vertices of $F$ is a subset of the set of vertices of $C$.
Then we write $F \cellleq C$ and call $C$ a supersimplex of $F$.
Accordingly, we write $F \cellless C$ for $F \cellleq C$ with $F \neq C$,
and $F \ncellleq C$ when $F \cellleq C$ does not hold.

We call a set of simplices $\calT$ a \emph{simplicial complex} provided that
for each simplex $C \in \calT$ all subsimplices of $C$ are included in $\calT$
and that any non-empty intersection of two simplices $C, C' \in \calT$ is a subsimplex of both $C$ and $C'$.

Given a simplicial complex $\calT$, we write $\calT^m = \{ C \in \calT \st \dim C = m \}$.
We say that $\calT$ is $p$-dimensional provided that
for each $S \in \calT$ there exists $C \in \calT^p$ with $S \cellleq C$.
An $m$-dimensional \emph{simplicial subcomplex} $\calU$ of $\calT$ is a subset of $\calT$ that is an $m$-dimensional simplicial complex by itself.
We write $\calT^{[m]}$ for the largest $m$-dimensional simplicial subcomplex of $\calT$,
which is called the \emph{$m$-skeleton} of $\calT$.
\\

Simplices are orientable compact manifolds with corners \cite[Chapter 10]{LeeSmooth}.
We henceforth assume that all simplices in $\calT$ are oriented,
and that this orientation is the Euclidean orientation for $n$-dimensional simplices.
Any oriented $m$-simplex $C$ induces an orientation on any $(m-1)$-dimensional subsimplex $F$,
and we then set $o(F,C)=1$ if either both orientations on $F$ coincide, or $o(F,C) = -1$ if those orientations differ.
\\

The space of \emph{simplicial $m$-chains} $\calC_{m}(\calT)$ is the real vector space generated by the set $\calT^{m}$.
It is easy to verify that the \emph{simplicial boundary operator} $\partial_{m} : \calC_m(\calT) \rightarrow \calC_{m-1}(\calT)$,
which is defined by
\begin{align}
 \partial_{m} C = \sum_{ \substack{ F \cellless C \\ F \in \calT^{m-1} } } o(F,C) F,
 \quad
 C \in \calT^{m},
\end{align}
satisfies $\partial_{m-1} \partial_m = 0$.
Thus we have a differential complex, 
\begin{align}
 \label{seq:prep_simplicial_calT}
 \begin{CD}
  0                      @>{              }>>
  \calC_{p  }(\calT)     @>{\partial_{p  }}>>
  \dots                  @>{\partial_{  1}}>>
  \calC_{  0}(\calT)     @>{              }>> 0,
 \end{CD}
\end{align}
called the \emph{simplicial chain complex} of $\calT$.

If $\calU$ is a simplicial subcomplex of $\calT$, then $\calC_m(\calU)$ is a subspace of $\calC_m(\calT)$,
and the simplicial boundary operator $\partial_{m} : \calC_m(\calU) \rightarrow \calC_{m-1}(\calU)$ is the restriction 
of the simplicial boundary operator $\partial_{m} : \calC_m(\calT) \rightarrow \calC_{m-1}(\calT)$.
This means that the simplicial chain complex of $\calU$,
\begin{align}
 \label{seq:prep_simplicial_calU}
 \begin{CD}
  0                      @>{              }>>
  \calC_{p  }(\calU)     @>{\partial_{p  }}>>
  \dots                  @>{\partial_{  1}}>>
  \calC_{  0}(\calU)     @>{              }>> 0,
 \end{CD}
\end{align}
is a differential subcomplex of \eqref{seq:prep_simplicial_calT}.

We are interested in the simplicial chain complex of $\calT$ relative to $\calU$.
In order to define that notion, let us write
\begin{align*}
 \calC_m(\calT,\calU) := \calC_m(\calT) / \calC_m(\calU)
\end{align*}
for the quotient space, and note that the equivalence classes
of the simplices $\calT^m \setminus \calU^m$ constitute a basis of $\calC_m(\calT,\calU)$.
We may therefore tacitly identify $\calC_m(\calT,\calU)$ with the vector space
generated by the set $\calT^{m} \setminus \calU^{m}$.
Now we observe for $C, C' \in \calC_m(\calT)$ with $C - C' \in \calC_m(\calU)$ that
\begin{align*}
   \partial_{m} C - \partial_{m} C'
 = \partial_{m} ( C - C' ) \in \calC_{m-1}(\calU).
\end{align*}
So the simplicial chain complexes \eqref{seq:prep_simplicial_calT} and \eqref{seq:prep_simplicial_calU}
induce another differential complex,
\begin{align}
 \label{seq:prep_simplicial_relative}
 \begin{CD}
  0                        @>{              }>>
  \calC_{p  }(\calT,\calU) @>{\partial_{p  }}>>
  \dots                    @>{\partial_{  1}}>>
  \calC_{  0}(\calT,\calU) @>{              }>> 0,
 \end{CD}
\end{align}
called the \emph{simplicial chain complex of $\calT$ relative to $\calU$}.
Up to elements of $\calU$, the differential of that complex is fully described by the equation
\begin{align}
 \partial_{m} C = \sum_{ \substack{ F \cellless C \\ F \in \calT^{m-1} \setminus \calU^{m-1} } } o(F,C) F,
 \quad
 C \in \calT^m \setminus \calU^m.
\end{align}
Note that \eqref{seq:prep_simplicial_relative} agrees with \eqref{seq:prep_simplicial_calT}
in the special case $\calU = \emptyset$.
The \emph{simplicial homology spaces} $\calH_m(\calT,\calU)$ 
of $\calT$ relative to $\calU$ are defined as the quotient spaces
\begin{align}
 \calH_m(\calT,\calU) &:=
  \dfrac{
  \ker \Big( \partial_{m} : \calC_{m}(\calT,\calU) \rightarrow \calC_{m-1}(\calT,\calU) \Big)
  }{
  \rng \Big(\partial_{m+1} : \calC_{m+1}(\calT,\calU) \rightarrow \calC_{m}(\calT,\calU) \Big)
  }.
\end{align}
Their dimensions are of general interest.
We call $b_m(\calT,\calU) := \dim \calH_m(\calT,\calU)$
the \emph{$m$-th simplicial Betti number} of $\calT$ relative to $\calU$,
and we call $b_m(\calT) := b_m(\calT,\emptyset)$ the \emph{$m$-th absolute simplicial Betti number} of $\calT$.
\\

Let $M$ be a topological manifold with boundary, embedded in $\bbR^n$,
and let $\Gamma$ be a topological submanifold of its boundary manifold $\partial M$.
We say a simplicial complex $\calT$ triangulates a topological manifold with boundary if that manifold is the union of all simplices in $\calT$.
The \emph{$m$-th topological Betti number} $b_m(M,\Gamma)$ is the dimension of the $m$-th singular homology group of $M$ relative to $\Gamma$;
we refer to \cite[Chapter 4, Section 4]{spanier1995algebraic} for the details.
In case $\Gamma = \emptyset$ we call $b_m(M) := b_m(M,\emptyset)$ the \emph{$m$-th absolute topological Betti number} of $M$.
The following canonical result relates the simplicial and the topological Betti numbers:

\begin{theorem}[{\cite[Chapter 4, Section 6, Theorem 8]{spanier1995algebraic}}]
 \label{ch_simplicial_thm_1}
 Let $M$ and $\Gamma$ be as in the previous paragraph.
 Let $\calT$ be a simplicial complex triangulating $M$,
 and $\calU$ be a simplicial subcomplex of $\calT$ that triangulates $\Gamma$.
 Then we have
 \begin{align*}
  b_m(\calT,\calU) = b_m(M,\Gamma)
 \end{align*}
 for all $m \in \bbN_{0}$. \qed
\end{theorem}

\begin{remark}
 This result implies that the topological Betti numbers can be calculated from the combinatorial structure of any triangulation.
 In certain applications they also coincide with the dimensions of the solution spaces of certain homogeneous partial differential equations over $M$;
 see Subsection \ref{subsect:preparation_derhamcomplex} below.
\end{remark}

\begin{example}
 \label{ex:betti_example}
 The following topological Betti numbers are of frequent interest.
 All Betti numbers $b_m(B^p)$ of the $p$-ball $B^p$ vanish except for $b_0(B^p) = 1$.
 All Betti numbers $b_m(S^p)$ of the $p$-sphere $S^p$ vanish except for $b_p(S^p) = b_0(S^p) = 1$.
 All Betti numbers $b_m(B^p,\partial B^p)$ of the $p$-ball relative to its boundary vanish except for $b_p(B^p,\partial B^p) = 1$.
 If $D^{p-1} \subsetneq \partial B^p$ is homeomorphic to $B^{p-1}$, 
 then all Betti numbers $b_m(B^p,D^{p-1})$ of the $p$-ball relative to a disk on the boundary vanish.
 \demo
\end{example}

\subsection{Differential forms on domains}
\label{subsect:preparation_differentialforms}

We review basic notions of differential forms on Riemannian manifolds with boundary,
and some aspects of their $L^2$ theory.
We generally refer to the textbooks \cite{gallot2004riemannian} and \cite{LeeSmooth} for general background on differential forms,
and to \cite{GMM} for more on their $L^2$ theory.
\\

Consider an $m$-dimensional open smooth manifold $M$ embedded in $\bbR^n$
such that its closure $\overline M$ is a topological manifold with boundary embedded in $\bbR^n$.
Examples include (the interiors of) simplices, polyhedrally bounded domains, and Lipschitz domains,
and moreover $\bbR^n$ itself.

Let $C^\infty\Lambda^k(M)$ be the space of smooth differential forms on $M$,
and let $C^\infty\Lambda^k(\overline M) \subseteq C^\infty\Lambda^k(M)$ be the image
of the pullback of $C^\infty\Lambda^k(\bbR^n)$ into $C^\infty\Lambda^k(M)$.
We let $C_c^\infty\Lambda^k(M)$ be the subspace of $C^\infty\Lambda^k(M)$
whose elements have compact support in $M$.
We know the \emph{exterior derivative}
\begin{align}
 \cartan_{M}^k &:
 C^\infty\Lambda^{k}(M)
 \rightarrow
 C^\infty\Lambda^{k+1}(M)
\end{align}
of a $k$-form, and and the \emph{$\wedge$-product}
\begin{align}
 \wedge &: 
 C^\infty\Lambda^{k}(M) \times C^\infty\Lambda^{l}(M)
 \rightarrow
 C^\infty\Lambda^{k+l}(M)
\end{align}
between a $k$-form and an $l$-form.
It is known for $\omega \in C^\infty\Lambda^k(M)$ and $\eta \in C^\infty\Lambda^l(M)$ that
\begin{align}
 \label{eq:wedge_properties}
 \omega \wedge \eta = (-1)^{kl} \eta \wedge \omega,
 \quad
 \cartan_{M}^{k+l}( \omega \wedge \eta ) = \cartan_{M}^{k}\omega \wedge \eta + (-1)^{k} \omega \wedge \cartan_{M}^{l} \eta.
\end{align}
The exterior derivative is linear, and it maps $C^\infty\Lambda^k(\overline M)$ into $C^\infty\Lambda^{k+1}(\overline M)$.
The $\wedge$-product is bilinear, and it maps $C^\infty\Lambda^{k}(\overline M) \times C^\infty\Lambda^{l}(\overline M)$ into $C^\infty\Lambda^{k+l}(\overline M)$.
\\

In the case that $M$ is oriented, the integral of $m$-forms over $M$ is well-defined,
and any non-vanishing $m$-form is either positively or negatively oriented.
If $M$ is moreover equipped with a Riemannian metric $g$,
then there exists a unique positively oriented normalized $m$-form $\vol_{M}$ over $M$,
called the \emph{volume form over $M$}.
The \emph{Hodge star operator} $\star_{M} : C^\infty\Lambda^{k}(\overline M) \rightarrow C^\infty\Lambda^{m-k}(\overline M)$
is uniquely defined by 
\begin{align}
 \label{eq:hodge_star_equation}
 \omega \wedge \star_{M} \eta
 =
 g(\omega,\eta) \vol_M
 ,
 \quad
 \omega, \eta \in C^\infty\Lambda^{k}(M).
\end{align}
For $\omega, \eta \in C^\infty\Lambda^{k}(\overline M)$ the Hodge star furthermore satisfies
\begin{align}
 \label{eq:hodge_star_fuss}
 \omega \wedge \star_{M}\eta = \eta \wedge \star_{M}\omega,
 \quad
 \star_{M}\star_{M} \omega = (-1)^{k(m-k)} \omega.
\end{align}
The \emph{exterior codifferential} is defined as
\begin{align}
 \label{eq:definition_codifferential}
 \delta_{M}^k : C^\infty\Lambda^{k}(M) \rightarrow C^\infty\Lambda^{k-1}(M),
 \quad
 \omega \mapsto (-1)^{m(k+1)+1} \star_{M} \cartan_{M}^{m-k} \star_{M} \omega.
\end{align}
Note that
\begin{align}
 \delta_{M}^k = (-1)^{k} \star_{M}^{\inv} \cartan_{M}^{m-k} \star_{M}.
\end{align}
If $\omega \in C^{\infty}\Lambda^k(M)$ and $\eta \in C^{\infty}\Lambda^{k+1}(M)$,
and at least one of these has support compact in $M$,
then the integration by parts formula
\begin{align}
 \label{eq:integration_by_parts_smooth}
 \int_M \cartan^{k}_{M} \omega \wedge \star_M \eta
 =
 \int_M \omega \wedge \delta^{k+1}_{M} \star_M \eta
\end{align}
holds. Given the Riemannian structure, we define the \emph{$L^2$ scalar product}
\begin{align}
 \langle \omega, \eta \rangle_{L^2\Lambda^k(M)}^2 
 &=
 \int_{\Omega} g(\omega, \eta ) \vol_M,
 \quad
 \omega, \eta \in C^{\infty}\Lambda^k(M)
 .
\end{align}
The completion of $C^{\infty}\Lambda^k(M)$ with respect to this scalar product
is the Hilbert space $L^{2}\Lambda^k(M)$ of square-integrable differential $k$-forms over $M$.
The Hodge star extends to a mapping $\star_{M} : L^2\Lambda^k(M) \rightarrow L^2\Lambda^{m-k}(M)$.
\\

Suppose that $\calT$ is a simplicial complex in $\bbR^{n}$.
We recall that we assume all simplices to have a fixed orientation.
We now furthermore suppose that $\bbR^n$ has a Riemannian metric,
which induces a Riemannian structure on every $C \in \calT$.
The notion of manifolds with corners \cite[Chapter 10]{LeeSmooth}
enables us to define traces along the boundary of simplices.

For all $C \in \calT$ and $F \cellleq C$ we have well-defined \emph{tangential trace operators}
\begin{align}
 \label{eq:tangential_trace_local}
 \trace^{k}_{C,F}       &: C^\infty\Lambda^k(C) \rightarrow C^\infty\Lambda^k(F),
\end{align}
They can be defined via the pullback of the inclusion of manifolds with corners.
We also find use for \emph{normal trace operators}, for which no canonical notation seems to exist. We write:
\begin{align}
 \label{eq:normal_trace_local}
 \normaltrace^{k}_{C,F}       &: C^\infty\Lambda^k(C) \rightarrow C^\infty\Lambda^{\dim F - \dim C + k}(F),
 \quad \omega \mapsto \star_{F}^\inv \trace^{\dim C-k}_{C,F} \star_{C} \omega.
\end{align}
These operators satisfy
\begin{align}
 \label{eq:trace_cartan_commute}
 \trace^{k+1}_{C,F}   \cartan_{C}^k &= \cartan^{k}_{F} \trace^{k}_{C,F},
 \\
 \label{eq:normaltrace_delta_commute}
 \normaltrace^{k-1}_{C,F} \cocartan_{C}^k &= (-1)^{\dim C - \dim F} \cocartan^{\dim F - \dim C + k}_{F} \normaltrace^{k}_{C,F}.
\end{align}
An important appearance of these operators is a variant of Stokes' theorem,
\begin{align}
 \label{eq:stokes_theorem}
 \begin{split}
 &\int_C \cartan^{k} \omega \wedge \star_C \eta
 -
 \int_C \omega \wedge \star_C \delta_{C}^{k+1} \eta
 \\=&
 \sum_{ \substack{ F \cellless C \\ F \in \calT^{\dim C - 1} } }
 o(F,C)
 \int_F \trace^{k}_{C,F} \omega \wedge \star_F \normaltrace^{k+1}_{C,F} \eta,
 \end{split}
\end{align}
for $\omega \in C^{\infty}\Lambda^{k}(C)$ and $\eta \in C^{\infty}\Lambda^{k+1}(C)$.
This generalizes the integration by parts formula \eqref{eq:integration_by_parts_smooth} above; 
see also \cite[Equation (0.2)]{dodziuk1976finite}.

Suppose that $\Omega$ is a polyhedral Lipschitz domain
carrying the Euclidean orientation and the Riemannian structure inherited from $\bbR^n$,
and that $\overline\Omega$ is triangulated by a finite simplicial complex $\calT$.
The tangential trace $\trace^{k}_{C} : C^\infty\Lambda^k(\overline\Omega) \rightarrow C^\infty\Lambda^k(C)$ for $C \in \calT$ is well-defined,
and we define the normal trace operator as
\begin{align}
\label{eq:normal_trace_global}
\normaltrace^{k}_{C}         &: C^\infty\Lambda^k(\overline\Omega) \rightarrow C^\infty\Lambda^{\dim C - n + k}(C),
\quad \omega \mapsto \star_{C}^\inv \trace^{n-k}_{C} \star_{\Omega} \omega
\end{align}
for notational convenience in the sequel.

\subsection{Hilbert complexes}
\label{subsect:preparation_hilbertcomplexes}

We review basic notions of Hilbert complexes,
which can be found in \cite{AFW2}.
We refer to \cite{bruening1992hilbert} for further background.
Our core results pertain to a class of merely finite-dimensional Hilbert complexes,
for which the following definitions simplify considerably,
but we frequently revisit an infinite-dimensional example application from
the theory of partial differential equations on manifolds,
which is described in the next subsection below.
\\

A \emph{Hilbert complex} is a sequence of real Hilbert spaces $X^{i}$, typically indexed over non-negative integers,
together with a sequence of closed densely-defined linear mappings $d^{i} : \dom(d^{i}) \subseteq X^{i} \rightarrow X^{i+1}$
which satisfy $\rng d^{i-1} \subseteq \ker d^{i}$.
\begin{align*}
 \begin{CD}
  0 @>>>
  X^{0  } @>{d^{0  }}>>
  X^{1  } @>{d^{1  }}>>
  \dots%
 \end{CD}
\end{align*}
Then the adjoint operators $d_i^\ast : \dom(d_i^\ast) \subseteq X^{i+1} \rightarrow X^i$
are densely-defined and closed as well.
We have the \emph{adjoint Hilbert complex}:
\begin{align*}
 \begin{CD}
  0 @<<<
  X^{0  } @<{d_{0  }^\ast}<<
  X^{1  } @<{d_{1  }^\ast}<<
  \dots%
 \end{CD}
\end{align*}
We assume that the operators $d^i$ have closed range,
and then $d_i^\ast$ have closed range, too.
Under this assumption, the space $\frakH^i = \ker d^i \cap \ker d_{i-1}^\ast$,
which we call \emph{$i$-th harmonic space}, satisfies
\begin{align}
 \frakH^i
 =
 \ker d^i \cap (\rng d^{i-1})^\perp
 =
 \ker d^\ast_{i-1} \cap (\rng d^\ast_{i})^\perp
 .%
\end{align}
We have the identities
\begin{align}
 \rng d^{i} = (\ker d_i^\ast)^\perp,
 \quad
 \rng d^\ast_{i} = (\ker d^i)^\perp.
\end{align}
The orthogonal decomposition
\begin{align}
 X^i = \rng d^{i-1} \oplus \rng d^\ast_{i} \oplus \frakH^i
\end{align}
is known as \emph{abstract Hodge decomposition} of $X^i$.
The $i$-th \emph{Hodge Laplacian} associated to a Hilbert complex
is the unbounded operator
\begin{align}
 \Delta_{i} := d^\ast_{i} d^{i} + d^{i-1} d^\ast_{i-1}.
\end{align}
It can be shown that $\Delta_i$ is closed, densely-defined, and has closed range.
In particular, the operator is self-adjoint, and 
\begin{align}
 \frakH^{i} = \ker \Delta_i = \left( \rng \Delta_i \right)^{\perp}.
\end{align}
The $i$-th Hodge Laplace problem associated to the Hilbert complex is then
to find $u \in \dom(\Delta_i)$ and $p \in \frakH^{i}$ such that
\begin{align}
 \Delta_i u = f - p, \quad u \perp \frakH^{i},
\end{align}
for given $f \in X^{i}$.
We refer to \cite{AFW2} for more on the variational and approximation theory of the Hodge Laplacian.

\subsection{The $L^2$ de Rham complex with partial boundary conditions}
\label{subsect:preparation_derhamcomplex}

Throughout this contribution, we develop an example application as a sideline,
where we demonstrate the general theory.
In this subsection, we introduce this basic example:
the $L^2$ de Rham complex with partial boundary conditions over a polyhedral Lipschitz domain.
The $L^2$ de Rham complex over $\Omega$ without boundary conditions is a prototypical example of a Hilbert complex.
We consider the general case of partial boundary conditions on the background of \cite{GMM},
which provides a theoretical setting for the Hodge Laplace equation with mixed boundary conditions.
\\

Let $\Omega$ be a bounded polyhedral Lipschitz domain.
The boundary $\partial\Omega$ is an $(n-1)$-dimensional Lipschitz manifold without boundary.
We assume that $\partial\Omega$ is the essentially disjoint union
of two $(n-1)$-dimensional Lipschitz manifolds $\Gamma_T$ and $\Gamma_N$ with boundary:
we call $\Gamma_T$ the tangential boundary part,
and we call $\Gamma_N$ the normal boundary part.

We furthermore assume that $\Omega$ is a triangulated by a finite simplicial complex $\calT$,
and that we have a simplicial subcomplex $\calU$ of $\calT$ that triangulates $\Gamma_N$.
Note that this implies that another simplicial subcomplex $\calV$ of $\calT$ triangulates $\Gamma_T$.
\\

We introduce spaces of smooth differential $k$-forms over $\Omega$ that satisfy partial boundary conditions.
The forms in the space $C_T^\infty\Lambda^k(\overline\Omega)$ satisfy tangential boundary conditions
along the tangential boundary part,
\begin{align*}
 C_T^\infty\Lambda^k(\overline\Omega)
 :=
 \left\{ 
  \omega \in C^\infty\Lambda^k(\overline\Omega) \st \forall F \in \calT, F \subseteq \Gamma_T : \trace^k_{F} \omega = 0
 \right\},
\end{align*}
and the forms in the space $C_N^\infty\Lambda^k(\overline\Omega)$ satisfy normal boundary conditions
along the normal boundary part,
\begin{align*}
 C_N^\infty\Lambda^k(\overline\Omega)
 &:=
 \left\{ 
  \omega \in C^\infty\Lambda^k(\overline\Omega) \st \forall F \in \calT, F \subseteq \Gamma_N : \normaltrace^{k}_{F} \omega = 0
 \right\}.
\end{align*}
Note that the integration by parts formula \eqref{eq:integration_by_parts_smooth} generalizes:
\begin{align*}
 \int_\Omega \cartan^{k}_{\Omega} \omega \wedge \star_\Omega \eta
 =
 \int_\Omega \omega \wedge \delta^{k+1}_{\Omega} \star_\Omega \eta,
 \quad
 \omega \in C_T^\infty\Lambda^k(\overline\Omega), \quad \eta \in C_N^\infty\Lambda^{k+1}(\overline\Omega).
\end{align*}
Defining $L^2$ differential forms that have distributional exterior derivative in $L^2$
and satisfy partial boundary conditions is far from trivial.
Given $\omega \in L^{2}\Lambda^k(\Omega)$, we write $\omega \in H\Lambda^{k}(\Omega)$
if there exists $\xi \in L^{2}\Lambda^{k+1}(\Omega)$ such that
\begin{align*}
 \left\langle \cartan^{k}_{\Omega} \omega, \eta \right\rangle_{L^2\Lambda^{k+1}(\Omega)}
 =
 \left\langle \xi, \delta^{k+1}_{\Omega} \eta \right\rangle_{L^2\Lambda^{k}(\Omega)},
 \quad
 \eta \in C_c^{\infty}\Lambda^{k+1}(\Omega),
\end{align*}
and in that case we define $\cartan^{k}_{\Omega} \omega = \xi$.
Similarly, we write $\omega \in H^{\star}\Lambda^{k}(\Omega)$
if there exists $\xi \in L^{2}\Lambda^{k-1}(\Omega)$ such that
\begin{align*}
 \left\langle \delta^{k}_{\Omega} \omega, \eta \right\rangle_{L^2\Lambda^{k-1}(\Omega)}
 =
 \left\langle \xi, \cartan^{k-1}_{\Omega} \eta \right\rangle_{L^2\Lambda^{k}(\Omega)}
 \quad
 \eta \in C_c^{\infty}\Lambda^{k-1}(\Omega),
\end{align*}
and in that case we define $\delta^{k}_{\Omega} \omega = \xi$.
It can be shown \cite{AFW1} that this defines closed densely-defined operators
\begin{align*}
 \cartan^{k} &: H        \Lambda^{k}(\Omega) \subseteq L^{2}\Lambda^{k}(\Omega) \rightarrow L^{2}\Lambda^{k+1}(\Omega),
 \\
 \delta^{k}  &: H^{\star}\Lambda^{k}(\Omega) \subseteq L^{2}\Lambda^{k}(\Omega) \rightarrow L^{2}\Lambda^{k-1}(\Omega),
\end{align*}
which have closed range.
In particular $H\Lambda^{k}(\Omega)$ and $H^{\star}\Lambda^{k}(\Omega)$ are Hilbert spaces.

Next we introduce partial boundary conditions.
For $\omega \in L^{2}\Lambda^k(\Omega)$ we let $\tilde\omega \in L^{2}\Lambda^{k}(\bbR^n)$
denote the trivial extension outside of $\Omega$.
We recall that for every $x \in \partial\Omega$ there exists a sufficiently small open ball $B_x \ni x$
such that $\partial\Omega$ splits $B_x$ into exactly two simply connected components.
We define $H_T\Lambda^k(\Omega)$ as the linear subspace of $H\Lambda^k(\Omega)$
such that $\omega \in H_T\Lambda^k(\Omega)$ if and only if
for all $x \in \Gamma_T$, the ball $B_x$ chosen so small that $B_x \cap \partial\Omega \subset \Gamma_T$,
and all $\phi \in C_c^\infty\Lambda^{k+1}(B_x)$ we have
\begin{align*}
 \int_{B_x} \widetilde{\cartan^{k} \omega} \wedge \star \phi
 =
 \int_{B_x} \widetilde\omega \wedge \star \delta^{k+1}\phi
 .
\end{align*}
Analogously, we define $H_N^\star\Lambda^k(\Omega)$ as the linear subspace of $H^\star\Lambda^k(\Omega)$
such that $\omega \in H_N^\star\Lambda^k(\Omega)$ if and only if
for all $x \in \Gamma_N$, the ball $B_x$ chosen so small that $B_x \cap \partial\Omega \subset \Gamma_N$,
and all $\phi \in C_c^\infty\Lambda^{k-1}(B_x)$ we have
\begin{align*}
 \int_{B_x} \widetilde{\delta^{k} \omega} \wedge \star \phi
 =
 \int_{B_x} \widetilde\omega \wedge \star \cartan^{k-1}\phi
 .
\end{align*}
It can easily be seen that $H_T\Lambda^k(\Omega)$ is a closed subspace of $H_T\Lambda^k(\Omega)$,
and that $H_N^\star\Lambda^k(\Omega)$ is a closed subspace of $H^\star\Lambda^k(\Omega)$.
By \cite[Theorem 4.4]{GMM} we have mutually adjoint closed densely-defined linear operators
\begin{align*}
 \cartan_T^{k} &: H_T        \Lambda^{k}(\Omega) \subseteq L^{2}\Lambda^{k}(\Omega) \rightarrow L^{2}\Lambda^{k+1}(\Omega),
 \\
 \delta_N^{k}  &: H_N^{\star}\Lambda^{k}(\Omega) \subseteq L^{2}\Lambda^{k}(\Omega) \rightarrow L^{2}\Lambda^{k-1}(\Omega).
\end{align*}
Moreover, the subspaces
\begin{align*}
 \cartan_T^k(H_T\Lambda^k) \subseteq H_T\Lambda^{k+1},
 \quad
 \delta_N^k(H^\ast_N\Lambda^k) \subseteq H^\ast_N\Lambda^{k-1}
\end{align*}
are closed \cite[Theorem 4.3]{GMM},
and in each of $H\Lambda^k(\Omega)$, $H_T\Lambda^k(\Omega)$, $H^\star\Lambda^k(\Omega)$, and $H_N^\star\Lambda^k(\Omega)$,
the respective intersection with $C^{\infty}\Lambda^{k}(\Omega)$ is dense \cite[Proposition 3.1]{jakab2009regularity}.
We conclude that we have mutually adjoint closed Hilbert complexes:
\begin{align}
\label{seq:l2_de_rham}
 \begin{CD} %
  0                                     @>>>
  H_T\Lambda^{0} \subseteq L^2\Lambda^0 @>\cartan_T^0>>
  \dots                                 @>\cartan_T^{n-1}>>
  H_T\Lambda^{n} \subseteq L^2\Lambda^n @>>>
  0
  \\ 
  0                                          @<<<
  H^\ast_N\Lambda^{0} \subseteq L^2\Lambda^0 @<\delta_N^1<<
  \dots                                      @<\delta_N^n<<
  H^\ast_N\Lambda^{n} \subseteq L^2\Lambda^n @<<<
  0 
 \end{CD}
\end{align}
The spaces of harmonic forms
\begin{align}
 \frakH^k(\Omega,\Gamma_T,\Gamma_N) = \ker \cartan_T^k \cap \ker \delta_N^k
\end{align}
of these two complexes are finite-dimensional \cite[Theorem 4.3]{GMM}, and
\begin{align}
\begin{split}
 b_k(\Omega,\Gamma_T) 
 &=
 \dim \frakH^k(\Omega,\Gamma_T,\Gamma_N)
 \\&=
 \dim \frakH^{n-k}(\Omega,\Gamma_N,\Gamma_T)
 =
 b_{n-k}(\Omega,\Gamma_N)
\end{split}
\end{align}
holds \cite[Theorem 5.3]{GMM}.
\\

A major motivation to study the $L^2$ de Rham complex with partial boundary conditions
is the Hodge Laplace equation over $k$-forms with mixed boundary conditions.
This is the $k$-th Hodge Laplace problem associated to \eqref{seq:l2_de_rham},
where for given $f \in L^2\Lambda^{k}(\Omega)$ we search
$u \in L^2\Lambda^{k}(\Omega)$ and $p \in \frakH^k(\Omega,\Gamma_T,\Gamma_N)$ such that
\begin{align*}
 u \in \dom( \cartan_T^{k} ) \cap \dom ( \delta_N^{k} ),
 \quad
 \delta_N^{k}  u \in \dom( \cartan_T^{k-1} ),
 \quad
 \cartan_T^{k} u \in \dom( \cartan_N^{k+1} ),
 \\
 \left( \delta_N^{k+1} \cartan_T^{k} + \cartan_T^{k-1} \delta_N^{k} \right) u + p = f,
 \quad
 u \perp \frakH^k(\Omega,\Gamma_T,\Gamma_N).
\end{align*}
The boundary conditions are interpreted in the way outlined above.
The complexes of finite element differential forms in this contribution pertain to numerical methods for this problem.

\ifx\justbeingincluded\undefined
\end{document}
\fi

\ifx\justbeingincluded\undefined
\input{../global/header_common.tex}

\begin{document}

\fi

\section{Discrete distributional differential forms}
\label{sect:element_systems}

This section introduces spaces of discrete distributional differential forms and operators mapping between those spaces.
We assume that $\calT$ is a finite $n$-dimensional simplicial complex and that $\calU$ is a subcomplex of $\calT$. 
We do not require a priori in this section that $\calT$ triangulates any manifold,
although it does in our example application from Subsection \ref{subsect:preparation_derhamcomplex}.
The subcomplex $\calU$ serves to formalize boundary conditions imposed on spaces of discrete distributional differential forms.
However, the reader may assume $\calU = \emptyset$ in a first reading.
\\

We consider finite element de Rham complexes on the simplices of $\calT$,
using assumptions which are inspired from abstract frameworks 
in the theory of finite element differential forms \cite{AFWgeodecomp,StructPresDisc}.
We assume that for each $m$-dimensional simplex $C \in \calT$ we have a sequence
of finite-dimensional subspaces $\Lambda^{k}(C)$ of $C^{\infty}\Lambda^{k}(C)$,
such that a finite element de Rham complex
\begin{align}
 \label{eq:simplex_sequence}
 \begin{CD}
   0 @>>> \Lambda^0(C)
     @>\Difftot_C^{  0}>> \Lambda^{  1}(C)
     @>\Difftot_C^{  1}>> \dots
     @>\Difftot_C^{m-1}>> \Lambda^{m  }(C)
     @>>> 0
  \end{CD}
\end{align}
is constituted. It is furthermore assumed that the trace induces
a surjective mapping from $\Lambda^k(C)$ to $\Lambda^k(F)$, i.e.
\begin{align*}
 \trace^k_{C,F} \Lambda^k(C) = \Lambda^k(F),
\end{align*}
for any $F, C \in \calT$ with $F \cellleq C$.

\begin{example}
 \label{example:dddf_anke}
 For our example application from Subsection \ref{subsect:preparation_derhamcomplex},
 consider a finite element complex of Arnold-Falk-Winther-type on $\bbR^n$, 
 for example the complex of trimmed polynomial differential forms of degree $r \geq 1$.
 \begin{align*}
  \begin{CD}
   \dots @>\Difftot^{k-1}>> \calP_r^-\Lambda^k(\bbR^n) @>\Difftot^{k}>> \calP_r^-\Lambda^{k+1}(\bbR^n) @>\Difftot^{k+1}>> \dots
  \end{CD}
 \end{align*}
 We then set $\Lambda^k(C)$ as the trace of $\calP_r^-\Lambda^k(\bbR^n)$ on the simplex $C \in \calT$.
 This results in discrete de Rham complexes on each simplex,
 and the traces are surjective.
 So our assumptions hold in this setting.
 The other Arnold-Falk-Winther-type complexes can be treated analogously,
 and we refer to \cite{AFW1,AFWgeodecomp} for further background.
 A finite element complex of non-uniform polynomial degree is described in \cite{zaglmayr2006high}
 and satisfies the above assumptions as well.
 \demo
\end{example}

Towards a definition of \emph{discrete distributional differential form},
let us consider $C \in \calT^m$ and $\omega_C \in \Lambda^k(C)$.
We then call $\omega_C$ a discrete distributional differential form of degree $k-n+m$.
The following example motivates this terminology.

\begin{example}
 \label{example:dddf_pia}
 We continue with our example application from Subsection \ref{subsect:preparation_derhamcomplex}.
 Let $C \in \calT^m$ and $\omega_C \in \Lambda^{k-n+m}(C)$.
 Then $\omega_C$ is a discrete distributional differential form of degree $k$,
 and $\omega_C$ can be interpreted as a linear functional on $C_N^\infty\Lambda^{k}(\overline\Omega)$
 by
 \begin{align*}
  \langle \omega_C, \phi \rangle
  :=
  \int_C\omega_C \wedge \star_C \normaltrace^{k}_C \phi,
  \quad
  \phi \in C_N^\infty\Lambda^{k}(\overline\Omega).
 \end{align*}
 We define the distributional exterior derivative $\cartan_{\Omega}^{k} \omega_C$ of $\omega_C$
 as a functional on $C_N^\infty\Lambda^{k+1}(\overline\Omega)$ via
 \begin{align*}
  \langle \cartan_{\Omega}^{k} \omega_C, \Phi \rangle
  :=
  \langle \omega_C, \delta_{\Omega}^{k+1} \Phi \rangle,
  \quad
  \Phi \in C_N^\infty\Lambda^{k+1}(\overline\Omega).
 \end{align*}
 Let us derive a different form for this functional.
 For $\Phi \in C_{N}^\infty\Lambda^{k+1}(\overline\Omega)$ we find
 \begin{equation*}
 \begin{split}
  \langle \omega_C, \delta_{\Omega}^{k+1} \Phi \rangle
  =
  \int_C \omega_C \wedge \star_C \normaltrace^{k}_C \delta_{\Omega}^{k+1} \Phi
  &=
  (-1)^{n-m}
  \int_C \omega_C \wedge \star_C  \delta_{C}^{k+1} \normaltrace^{k+1}_C \Phi,
 \end{split}
 \end{equation*}
 using \eqref{eq:normaltrace_delta_commute}.
 Applying the integration by parts formula \eqref{eq:stokes_theorem}
 and using the elementary observation
 $\normaltrace^{k+m-n+1}_{C,F} \normaltrace^{k+1}_C = \normaltrace^{k+1}_F$,
 we see that
 \begin{equation*}
  \begin{split}
   (-1)^{n-m}
   \langle \omega_C, \delta_{\Omega}^{k+1} \Phi \rangle
   =
   \left\langle \cartan_{C}^{k-n+m} \omega, \Phi \right\rangle
   -
   \sum_{ \substack{ F \cellless C \\ F \in \calT^{m - 1} } }
   o(F,C) \left\langle \trace^{k-n-m}_{C,F}\omega, \Phi \right\rangle.
  \end{split}
 \end{equation*}
 From $\delta_{\Omega}^{k} \delta_{\Omega}^{k+1} = 0$, or through direct computation,
 we obtain $\cartan_{\Omega}^{k+1} \cartan_{\Omega}^{k} = 0$ for the distributional exterior derivative.
 \demo
\end{example}

With this example in mind, we want to define the distributional exterior derivative of a discrete distributional differential form
in purely discrete terms without reference to a space of test functions.
Furthermore, we have seen in the example that the distributional derivative of a discrete distributional differential form
is the sum of a piecewise exterior derivative, and the sum of signed traces onto lower dimensional simplices.
We identify these two operations as operators in their own right.

In order to find suitable spaces between which these operators can map,
let us consider the direct sum
\begin{align}
 \label{eq:dddf_definition_atomic}
 \Lambda_{-1}^k(\calT^m,\calU)
 :=
 \bigoplus_{ C \in \calT^m \setminus \calU^m } \Lambda^{k}(C),
\end{align}
which gathers the $k$-forms associated to $m$-simplices.
For $\omega \in \Lambda_{-1}^k(\calT^m,\calU)$ we then write $\omega_{C}$
for the component of $\omega$ in \eqref{eq:dddf_definition_atomic}
corresponding to $C \in \calT^m \setminus \calU^m$.

We introduce a family of operators
\begin{align}
 \Diffhor^m_k :
 \Lambda^{k  }_{-1}(\calT^m,\calU)
 \longrightarrow
 \Lambda^{k+1}_{-1}(\calT^m,\calU),
\end{align}
called \emph{horizontal differentials}, defined by
\begin{align}
 \Diffhor^m_k \omega
 :=
 \sum_{C \in \calT^m \setminus \calU^m} \Difftot^k_C \omega_C,
\end{align}
and a family of differential operators
\begin{align}
 \Diffver^m_k :
 \Lambda^{k}_{-1}(\calT^{m  },\calU)
 \longrightarrow
 \Lambda^{k}_{-1}(\calT^{m-1},\calU),
\end{align}
called \emph{vertical differentials}, defined by
\begin{align}
 \Diffver^m_k \omega
 :=
 \sum_{ C \in \calT^m \setminus \calU^m}
 \sum_{ \substack{ F \cellless C \\ F \in \calT^{m-1} \setminus \calU^{m-1} } }
 o(F,C) \trace^k_{C,F} \omega_C.
\end{align}
The families of operators $\Diffhor^m_k$ and $\Diffver^m_k$ correspond to the two terms
that we identified in Example \ref{example:dddf_pia} above.
Before we properly the define the discrete distributional exterior derivative,
we introduce two more families of spaces:
\begin{align}
 \label{eq:dddf_definition_primal}
 \Lambda_{-b}^k(\calT^m,\calU)
 &:=
 \bigoplus_{j=0}^{b-1} \Lambda^{k-j}_{-1}(\calT^{m-j},\calU),
 \quad
 0 \leq b \leq m+1,
 \\
 \label{eq:dddf_definition_dual}
 \Kappa^k_{-b}(\calT^m,\calU)
 &:=
 \bigoplus_{j=0}^{b-1} \Lambda^{k+j}_{-1}(\calT^{m+j},\calU),
 \quad
 0 \leq b \leq n-m+1.
\end{align}
We henceforth agree that any element of a space $\Lambda_{-b}^k(\calT^m,\calU)$ or $\Kappa_{-b}^k(\calT^m,\calU)$
can be called a \emph{discrete distributional differential form}.
Every $\omega \in \Lambda_{-b}^k(\calT^m,\calU)$ gives a linear functional on $C_N^{\infty}\Lambda^{k}(\overline\Omega)$,
and therefore we say that $\omega$ is a discrete distributional differential form of degree $k$.
Analogously every $\omega \in \Kappa^k_{-b}(\calT^m,\calU)$ gives a linear function on $C_N^{\infty}\Lambda^{n-m}(\overline\Omega)$
and is thus called a discrete distributional differential form of degree $n-m$.
We note furthermore that $\Kappa^k_{-1}(\calT^m,\calU) = \Lambda^k_{-1}(\calT^m,\calU)$.
As we see below, the family $\Lambda_{-b}^{k}(\calT^{m},\calU)$ generalizes the standard finite element spaces,
whereas the family $\Kappa_{-b}^k(\calT^{m},\calU)$ generalizes the spaces of simplicial chains.

We use the term \emph{discrete distributional exterior derivative} for the linear mappings
\begin{align}
 \Difftot^{k+n-m} &: \Lambda^k_{-b}(\calT^m,\calU) \rightarrow \Lambda^{k+1}_{-b-1}(\calT^m,\calU),
 \\
 \Difftot^{k+n-m} &: \Kappa^k_{-b}(\calT^m,\calU) \rightarrow \Kappa^{k}_{-b-1}(\calT^{m-1},\calU),
\end{align}
which are uniquely defined by
\begin{align}
 \label{eq:distributional_differential}
 \Difftot^{k} \omega = (-1)^{i} \Diffhor^{n-i}_{k-i} \omega - (-1)^{i}\Diffver^{n-i}_{k-i} \omega,
 \quad
 \omega \in \Lambda^{k-i}_{-1}(\calT^{n-i},\calU).
\end{align}
The differential properties
\begin{align}
 \label{eq:differential_properties}
 \Diffhor^{k+1}_m \Diffhor^k_m = 0,
 \quad
 \Diffver^k_{m-1} \Diffver^k_m = 0,
 \quad
 \Difftot^{k+1} \Difftot^k = 0
\end{align}
can be verified by direct computation.
\\

%

The kernels of $\Diffhor^m_k$ and $\Diffver^m_k$ are interesting in their own right.
We define
\begin{align}
 \label{eq:dddf_mathieu}
 \Lambda^k(\calT^m,\calU) 
 &:=
 \left\{
  \omega \in \Lambda^k_{-1}(\calT^m,\calU) \;\middle|\; \Diffver^m_k \omega = 0
 \right\},
 \\
 \label{eq:dddf_marina}
 \Kappa^k(\calT^m,\calU) 
 &:=
 \left\{
  \omega \in \Lambda^k_{-1}(\calT^m,\calU) \;\middle|\; \Diffhor^m_k \omega = 0
 \right\}.
\end{align}
We have well-defined operators
\begin{align}
 \cartan^{k+n-m} &: \Lambda^k(\calT^m,\calU) \rightarrow \Lambda^{k+1}(\calT^m,\calU),
 \\
 \cartan^{k+n-m} &: \Kappa^k(\calT^m,\calU) \rightarrow \Kappa^{k}(\calT^{m-1},\calU).
\end{align}
We may write
\begin{align*}
 \Lambda_{0}^k(\calT^m,\calU) = \Lambda^k(\calT^m,\calU),
 \quad
 \Kappa_{0}^k(\calT^m,\calU) = \Kappa^k(\calT^m,\calU)
\end{align*}
in order to unify the notation.
\\

The subcomplex $\calU$ plays only a minor technical role in our derivations.
In order to slim down the notation in the sequel, we therefore assume that $\calU$ is understood,
and we merely write
\begin{align*}
 \begin{array}{ll}
 \Lambda^k(\calT^m)
 \equiv
 \Lambda^k(\calT^m, \calU)
 ,
 &
 \Lambda_{-b}^k(\calT^m)
 \equiv
 \Lambda_{-b}^k(\calT^m, \calU)
 ,
 \\
 \Kappa^k(\calT^m)
 \equiv
 \Kappa^k(\calT^m, \calU)
 ,
 &
 \Kappa_{-b}^k(\calT^m)
 \equiv
 \Kappa_{-b}^k(\calT^m, \calU)
 \end{array}
\end{align*}
in the sequel.

\begin{example}
 Since these definitions are comparatively abstract,
 let us briefly anticipate the meaning of $\Lambda^k(\calT^n)$
 and $\Kappa^k(\calT^n)$ in our example application.
 
 The space $\Lambda^k(\calT) = \Lambda^k(\calT,\calU)$ corresponds
 to the standard finite element complexes with partial boundary conditions
 along the subcomplex $\calV$ which triangulates the tangential boundary $\Gamma_T$.
 To see this, note that the fact $\Diffver^n_k \omega = 0$ for $\omega \in \Lambda_0^k(\calT,\calU)$
 and Euclidean orientation of each $n$-simplex
 imply that along interior faces the tangential traces of neighboring $n$-simplices coincide,
 while along faces of $\calV$ the tangential traces must vanish.
 
 In our applicational setting, 
 the space $\Kappa^0(\calT^m)$ corresponds to $\calC_{m}(\calT,\calU)$.
 \demo
\end{example}

The distributional exterior derivative generalizes the classical exterior de\-ri\-va\-tive.
Since the families of operators $\Diffhor^k_m$, $\Diffver^k_m$, and $\Difftot^k$ satisfy the differential properties \eqref{eq:differential_properties},
we are motivated to consider differential complexes of discrete distributional differential forms.
%
But differential complexes are defined in purely algebraic terms.
Since Hilbert space structures are relevant in finite element theory additionally,
we utilize the theory of Hilbert complexes,
and assume that we have chosen a scalar product on the finite-dimensional spaces $\Lambda^k_{-1}(\calT^m)$.
The choice of the scalar product might depend on the respective application,
and determine the harmonic spaces of the Hilbert complexes but not their dimension.

\begin{example}
 As a possible choice of a scalar product on $\Lambda^k_{-1}(\calT^m)$,
 we consider the mesh-dependent scalar product
 \begin{align}
  \label{eq:dddf_scalar_product}
  \langle \omega, \eta \rangle
  :=
  \sum_{ C \in \calT^m - \calU^m } 
  h_C^{n-m} \langle \omega_C, \eta_C \rangle_{L^2\Lambda^{k}(C)},
  \quad
  \omega, \eta \in \Lambda^k_{-1}(\calT^m).
 \end{align}
 Here, $h_C$ denotes the diameter of $C$ if $\dim C \geq 1$, and, say, the average diameter of all adjacent edges if $C$ is a vertex.
 This scalar product appears in literature on residual error estimators \cite{Braess2007,BrSchoMax},
 but also in relation to discontinuous Galerkin methods \cite[Equation (2.2)]{arnold1982interior}.
 The weighting factors in \eqref{eq:dddf_scalar_product} are relevant for scaling arguments.
 \demo
\end{example}

\begin{remark}
 Our notion of discrete distributional differential form is similar but different from the notion of currents \cite{rham1984differentiable} introduced by de Rham.
 Currents in the sense of de Rham over an $n$-dimensional manifold are functionals on smooth differential $k$-forms,
 and they generalize differential $(n-k)$-forms via the pairing of $(n-k)$-forms with $k$-forms.
 Thus their definition only involves the oriented smooth structure on the manifold,
 whereas our notion of discrete distributional differential form is based on the $L^2$ pairing,
 which requires additionally a Riemannian structure on the manifold.
\end{remark}

\ifx\justbeingincluded\undefined
\end{document}
\fi



\ifx\justbeingincluded\undefined
\input{../global/header_common.tex}

\begin{document}

\fi

\section{Horizontal and Vertical Homology Theory}
\label{sect:homology_theory_theoretical}

In this section, Hilbert complexes with the differential operators $\Diffhor^m_k$ and $\Diffver^m_k$,
and a double complex are studied.
Under conditions that hold in applications,
we provide a new proof for the spaces of discrete harmonic $k$-forms of the standard finite element complex
to have dimension equal to the $k$-th Betti numbers of the triangulation.
Since we allow for partial boundary conditions,
we also close a gap in the literature on finite element differential forms.
\\

We continue to assume that we have a finite $n$-dimensional simplicial complex $\calT$,
a subcomplex $\calU$, and on each simplex $C \in \calT$
a finite-dimensional smooth de Rham sequence
such that the traces are surjective.
\\

Let us introduce more Hilbert complexes:
one with the differentials $\Diffhor^{m}_{k}$ for $m$ fixed,
and one with the differentials $\Diffver^{m}_{k}$ for $k$ fixed.
When $m$ is fixed, then we may consider the \emph{horizontal} Hilbert complex
\begin{align}
 \label{eq:dddf_horizontal_sequence}
 \begin{CD}
  0
  @>>>
  \Lambda^0_{-1}(\calT^m)
  @>{\Diffhor^{m}_{0}}>>
  \Lambda^{1}_{-1}(\calT^{m})
  @>{\Diffhor^{m}_{1}}>>
  \dots
  @>{\Diffhor^{m}_{m-1}}>>
  \Lambda^{m}_{-1}(\calT^{m})
  @>>>
  0,
 \end{CD}
\end{align}
and when $k$ is fixed, then we may consider the \emph{vertical} Hilbert complex
\begin{align}
 \label{eq:dddf_vertical_sequence}
 \begin{CD}
  0
  @>>>
  \Lambda^k_{-1}(\calT^{n  })
  @>{\Diffver^{n  }_{k}}>>
  \Lambda^k_{-1}(\calT^{n-1})
  @>{\Diffver^{n-1}_{k}}>>
  \dots
  @>{\Diffver^{k}_{k+1}}>>
  \Lambda^k_{-1}(\calT^{k})
  @>>>
  0.
 \end{CD}
\end{align}
We write $\frakH^k_\Diffhor(\calT^m)$ for the harmonic spaces of the Hilbert complex \eqref{eq:dddf_horizontal_sequence},
called \emph{horizontal harmonic spaces},
\begin{align}
 \frakH^k_\Diffhor(\calT^m) &:= \left\{ \omega \in \Lambda_{-1}^{k}(\calT^{m}) \st
 \omega \in \ker \Diffhor^{m}_{k}, \omega \perp \Diffhor^{m}_{k-1}\Lambda_{-1}^{k}(\calT^{m})
 \right\},
\end{align}
and we write $\frakH^k_\Diffver(\calT^m)$ for the harmonic spaces of the Hilbert complex \eqref{eq:dddf_vertical_sequence},
called \emph{vertical harmonic spaces}.
\begin{align}
 \frakH^k_\Diffver(\calT^m) &:= \left\{ \omega \in \Lambda_{-1}^{k}(\calT^{m}) \st
 \omega \in \ker \Diffver^{m}_{k}, \omega \perp \Diffhor^{m+1}_{k}\Lambda_{-1}^{k}(\calT^{m+1})
 \right\}.
\end{align}
We note that
\begin{align}
 \frakH^0_\Diffhor(\calT^m) = \Kappa^0(\calT^m),
 \quad
 \frakH^k_\Diffver(\calT^n) = \Lambda^k(\calT^n)
\end{align}
by definition.
In order to deduce more information about these complexes,
we make additional assumptions on the finite element spaces and the combinatorial properties of $\calT$,
to be described below.
\\


First, the horizontal complexes and the harmonic spaces $\frakH_\Diffhor^k(\calT^m)$ are considered. 
We recall that the complex
\begin{align}
 \label{eq:extended_horizontal_sequence}
 \begin{CD}
  0
  @>>>
  \Kappa^0(\calT^m)
  @>{}>>
  \Lambda^0_{-1}(\calT^m)
  @>{\Diffhor^m_{0  }}>>
  \dots
  @>{\Diffhor^m_{m-1}}>>
  \Lambda^m_{-1}(\calT^m)
  @>>>
  0
 \end{CD}
\end{align}
is
the direct sum of the simplex-wise complexes
\begin{align}
 \label{eq:simplex_sequence_extended}
 \begin{CD}
  0
  @>>>
  \ker \cartan_C^0
  @>{}>>
  \Lambda^{0  }(C)
  @>{\cartan_C^{0  }}>>
  \dots
  @>{\cartan_C^{m-1}}>>
  \Lambda^{m}(C)
  @>{}>>
  0
 \end{CD}
\end{align}
over $m$-simplices $C \in \calT^m \setminus \calU^{m}$.

\begin{definition}
We say that the \emph{local exactness condition} holds
if for each $C \in \calT \setminus \calU$, the sequence \eqref{eq:simplex_sequence_extended} is exact,
and if furthermore $\ker \cartan_{C}^0$ is spanned by the indicator function $1_C$ of $C$.
\end{definition}

This condition means that \eqref{eq:simplex_sequence} realizes the absolute cohomology on each simplex of $\calT \setminus \calU$.
It implies that $\frakH^k_\Diffhor(\calT^m)$ is trivial for $k \geq 1$,
and that $\frakH^0_\Diffhor(\calT^m)$ is spanned by the local indicator functions $1_C$, ${C \in \calT^m - \calU^m}$.
If the local exactness condition holds, then also
\begin{align*}
 \frakH_\Diffhor^0(\calT^m) \simeq \calC_m(\calT,\calU)
\end{align*}
by identifying each simplex with its local indicator function.

\begin{example}
 \label{ex:local_exactness_condition}
 We continue our example application.
 The finite element complexes of finite element exterior calculus \cite{AFW1} satisfy the local exactness condition.
 The finite element complexes of non-uniform polynomial degree developed in \cite{zaglmayr2006high} satisfy the local exactness condition as well.
 \demo
\end{example}


Next, we consider the vertical complexes and the vertical harmonic spaces $\frakH_\Diffver^k(\calT^m)$. 
We want to find conditions under which the sequence
\begin{align}
 \label{eq:trace_complex_extended}
 \begin{CD}
  0
  @>{}>>
  \Lambda^k(\calT)
  @>{}>>
  \Lambda_{-1}^k(\calT^n)
  @>{\Diffver^n_k}>>
  \dots
  @>{\Diffver^{k+1}_k}>>
  \Lambda_{-1}^k(\calT^k)
  @>>>
  0
 \end{CD}
\end{align}
is exact.
We show that this complex is the direct sum of complexes associated to local element patches,
after assuming a condition on the spaces $\Lambda^k(C)$.
Then we show exactness of those local sequences,
which requires another condition on the triangulation.
The exactness of \eqref{eq:trace_complex_extended} can then be concluded.
\\

In order to construct local vertical complexes associated to patches,
we use the notion of geometric decomposition of finite element spaces
which has been developed in \cite{AFW1} and \cite{AFWgeodecomp}.
For $C \in \calT^m$ let
\begin{align}
 \mathring\Lambda^{k}(C) := \left\{ \omega \in \Lambda^{k}(C) \st \forall F \cellless C : \trace_{C,F} \omega = 0 \right\}
\end{align}
denote the subspace of $\Lambda^k(F)$
whose members have vanishing traces on the boundary simplices of $F$.
\begin{definition}
We say the \emph{geometric decomposition condition} holds
if we have linear extension operators
\begin{align*}
 \ext^k_{F,C} : \mathring\Lambda^k(F) \rightarrow \Lambda^k(C),
\end{align*}
for $F, C \in \calT$ with $F \cellleq C$, such that
$\trace^k_{C,F} \ext_{F,C}^k = \Id_{\mathring\Lambda^k(F)}$ and
\begin{align*}
 \ext^{k}_{F,G} &= \trace^{k}_{C,G} \ext^{k}_{F,C},
 \quad G \in \calT, \quad F \cellleq G \cellleq C,
 \\
 0 &= \trace^{k}_{C,G} \ext^{k}_{F,C},
 \quad G \cellleq C, \quad F \ncellleq G,
\end{align*}
holds.
\end{definition}
Under this conditions one can derive the eponymous geometric decomposition of the distributional finite element spaces:
\begin{align}
 \label{eq:geometric_decomposition}
 \Lambda^k_{-1}(\calT^m)
 = 
 \bigoplus_{ C \in \calT^m \setminus \calU^m }
 \bigoplus_{ F \cellleq C }
 \ext^{k}_{F,C} \mathring\Lambda^{k}(F).
\end{align}
This is easy to see and follows also from a careful reading of \cite{AFWgeodecomp}. The authors make stronger assumptions in that publication,
but their derivation of \eqref{eq:geometric_decomposition} only requires the conditions that we assume in this work.
Extension operators that satisfy only our weaker assumptions have been also been used in Subsections 4.7 and 4.8 of \cite{AFW1}.
We refer to \cite[Proposition 2.2]{christiansen2013high} for a similar result.

\begin{example}
 \label{ex:compatible_extension_condition}
 We continue our example application.
 The geometric decomposition condition holds for the ansatz spaces of finite element exterior calculus;
 see Theorem 4.3, Theorem 7.3 and Theorem 8.3 of \cite{AFWgeodecomp}.
 It is also used implicitly for the complexes of non-uniform polynomial degree in \cite{zaglmayr2006high}.
 \demo
\end{example}

Assuming that the geometric decomposition assumption holds, we can now construct the local vertical complexes.
We define the spaces
\begin{align}
 \Kappa^m_k(F)
 :=
 \bigoplus_{ \substack{ C \in \calT^m \setminus \calU^m \\ F \cellleq C } }
 \ext^k_{F,C} \mathring\Lambda^{k}(F),
 \quad
 F \in \calT,
\end{align}
and note that $\Diffver^m_k \Kappa^m_k(F) \subseteq \Kappa^{m-1}_k(F)$. It is easy to verify:

\begin{lemma}
 Assume the geometric decomposition condition holds.
 Then the complex
 \begin{align}
  \begin{CD}
    0
    @>{}>>
    \Lambda^k(\calT)
    @>{}>>
    \Lambda_{-1}^k(\calT^n)
    @>{\Diffver^n_k}>>
    \dots
    @>{\Diffver^{k+1}_k}>>
    \Lambda_{-1}^k(\calT^k)
    @>>>
    0
  \end{CD}
 \end{align}
 is the direct sum of the complexes
 \begin{align*}
  \begin{CD}
    0
    @>{}>>
    \Kappa^n_k(F) \cap \ker {\Diffver^n_k}
    @>{}>>
    \Kappa^n_k(F)
    @>{\Diffver^n_k}>>
    \dots
    @>{\Diffver^{k+1}_k}>>
    \Kappa^k_k(F)
    @>>>
    0
  \end{CD}
 \end{align*}
 over all $F \in \calT$.
\end{lemma}

\begin{proof}
 Suppose that
 \begin{align*}
  \omega = \sum_{ C \in \calT^m \setminus \calU^m } \omega_C \in \Lambda_{-1}^{k}(\calT^m),
  \quad
  \omega_C \in \Lambda^k(C).
 \end{align*}
 By the geometric decomposition condition, we have
 \begin{align*}
  \omega_C = \sum_{ F \cellleq C } \omega_C^{(F)},
  \quad
  \omega_C^{(F)} \in \ext^k_{F,C} \mathring\Lambda^{k}(F).
 \end{align*}
 We rearrange this sum:
 \begin{align*}
  \omega = \sum_{ F \in \calT^{[m]} } \omega^{(F)},
  \quad
  \omega^{(F)}
  =
  \sum_{ \substack{ C \in \calT^m \setminus \calU^{m} \\ F \cellleq C } } \omega_C^{(F)}
  \in
  \Kappa^m_k(F).
 \end{align*}
 The $\omega^{(F)}$ provide the desired decomposition of $\Lambda_{-1}^{k}(\calT^m)$.
 Furthermore the geometric decomposition condition implies
  $\trace^{C,G} \omega_C^{(F)} \in \ext^k_{F,G} \mathring\Lambda^{k}(F)$
 for any $G \in \calT^{m-1}\setminus\calU^{m-1}$
 with $G \cellleq C$ and $F \cellleq G$. 
 This completes the proof.
\end{proof}

In order to prove the exactness of the local vertical sequences, 
we need to assume that the local combinatorial structure of $\calT$ and $\calU$ is sufficiently simple.
In order to verbalize a precise condition, we define for $F \in \calT$ the simplicial complexes
\begin{align}
 \begin{split}
 \label{eq:patch_complexes}
 \calM_F &:=
 \left\{ G \in \calT \st \exists C \in \calT^n : F \cellleq C \text{ and } G \cellleq C \right\},
 \\
 \calN_F &:= \left\{
  G \in \calM_F^{[n-1]} \;\middle|\; F \ncellleq G \text{ or } G \in \calU
 \right\}.
 \end{split}
\end{align}
The underlying idea is that $\calM_F$ triangulates the element patch around $F$,
and that $\calN_{F}$ triangulates the boundary of that patch, with modifications at $\calU$.
The condition that we assume is:
\begin{definition}
We say that $\calT$ satisfies the \emph{local patch condition} relative to $\calU$,
if for all $F \in \calT$ the complex
\begin{align}
 \label{eq:local_patch_sequence}
 \begin{CD}
  0 @>>>
  \calC_{n  }(\calM_F,\calN_F) @>{\partial_n}>>
  \dots                        @>{\partial_1}>>
  \calC_{0  }(\calM_F,\calN_F) @>>>
  0
 \end{CD}
\end{align}
has vanishing homology spaces at indices $n-1$, $\dots$, $0$.
\end{definition}
%
%

\begin{example}
 \label{ex:local_patch_condition}
 We continue with our example application,
 and show that the local patch condition holds there.
 Recall that we assume the boundary $\partial\Omega$
 to be decomposed into a normal boundary part $\Gamma_N$, triangulated by the subcomplex $\calU$,
 and a tangential boundary part $\Gamma_T$, triangulated by a subcomplex $\calV$.
 
 For $0 \leq d < \dim F,$ we observe that $\calM_F^d = \calN_F^d$ by definition.
 So in these cases we have $\calC_d(\calM_F,\calN_F) = 0$, and accordingly $b_d(\calM_F,\calN_F) = 0$.
 
 Now consider the case $d = \dim F$. 
 If $F \notin \calU$, then $\calM_F^{\dim F} \setminus \calN_F^{\dim F} = \{F\}$,
 and any simplex in $\calM_F$ of dimension $\dim F + 1$ has $F$ as a face;
 if instead $F \in \calU$, then $\calC_{\dim F}(\calM_F,\calN_F) = 0$.
 This means that \eqref{eq:local_patch_sequence} is always exact at $0 \leq d \leq \dim F$.
 
 Lastly, we treat the cases $\dim F < d < n$.
 Note that $\calM_F$ always triangulates a topological ball containing $F$.
 If $F$ is not a boundary simplex,
 then $\calN_F$ is just the boundary of that ball,
 and the exactness of \eqref{eq:local_patch_sequence} at $\dim F < d < n$ follows from Example \ref{ex:betti_example}.
 Note that $b_n(\calM_F,\calN_F) = 1$ corresponds to the contribution of $\mathring\Lambda^k(F)$ to the global finite element space $\Lambda^k(\calT,\calU)$.
 If instead $F$ is a boundary simplex with $F \notin \calV$,
 then the same argument applies.
 If $F$ is a boundary simplex with $F \in \calV$,
 then $\calN_F$ is a ball-shaped patch on the boundary of $\calM_F$,
 and $b_d(\calM_F,\calN_F)$ for $\dim F < d \leq n$.
 Note that $b_n(\calM_F,\calN_F) = 0$ corresponds to the vanishing trace
 of the global finite element space $\Lambda^k(\calT,\calU)$ on simplices of $\calV$ in that case.
 \demo
\end{example}

\begin{lemma}
 \label{res:local_patch_condition_skeleton}
 If the local patch condition holds for an $n$-dimensional simplicial complex $\calT$
 relative to an $(n-1)$-dimensional subcomplex $\calU$,
 then the local patch condition holds for $\calT^{[n-1]} \setminus \calU^{n-1}$ relative to $\calU^{[n-2]}$.
\end{lemma}

\begin{proof}
 Let $F \in \calT^{[n-1]} \setminus \calU^{n-1}$,
 let $\calM_F$ and $\calN_F$ as in \eqref{eq:patch_complexes},
 and let $\calM'_F$ and $\calN'_F$ be the corresponding subcomplexes
 in the formulation of the local patch condition of 
 $\calT^{[n-1]} \setminus \calU^{n-1}$ relative to $\calU^{[n-2]}$.
 We observe that $\calC_{m}(\calM_F,\calN_F)$, where $0 \leq m \leq n-1$, is spanned
 by the equivalence classes of the simplices $G \in \calT^{m} - \calU^{m}$
 with $F \cellleq G$.
 With similar reasoning, we observe that $\calC_{m}(\calM'_F,\calN'_F)$, where $0 \leq m \leq n-1$,
 is spanned by precisely the equivalence classes of the simplices.
 Thus we have $\calC_{m}(\calM_F,\calN_F) = \calC_{m}(\calM'_F,\calN'_F)$ for $0 \leq m \leq n-1$
 from which the statement follows.
\end{proof}

%

\begin{lemma}
 Assume the geometric decomposition condition
 and the local patch condition hold.
 Then the sequence
 \begin{align}
  \label{seq:sub}
  \begin{CD}
    0
    @>{}>>
    \Kappa^n_k(F) \cap \ker {\Diffver^n_k}
    @>{}>>
    \Kappa^n_k(F)
    @>{\Diffver^n_k}>>
    \dots
    @>{\Diffver^{k+1}_k}>>
    \Kappa^k_k(F)
    @>>>
    0
  \end{CD}
 \end{align}
 is exact for each $F \in \calT$.
\end{lemma}

\begin{proof}
 Let $F \in \calT$, and consider the differential complex
 that is derived from \eqref{eq:local_patch_sequence} by taking the tensor product $\mathring\Lambda^{k}(F)$
 at every index, i.e.
 \begin{align}
  \label{seq:dom}
  \begin{CD}
   \dots
   @>{ \partial_{m+1} \otimes \Id }>>
   \calC_{m  }(\calM_F,\calN_F) \otimes \mathring\Lambda^{k}(F)
   @>{ \partial_{m  } \otimes \Id }>>
   \calC_{m-1}(\calM_F,\calN_F) \otimes \mathring\Lambda^{k}(F)
   @>{ \partial_{m-1} \otimes \Id }>>
   \dots
  \end{CD}
 \end{align}
 The $m$-th homology space of \eqref{seq:dom}
 is isomorphic to $\calH_{m}(\calM_F,\calN_F) \otimes \mathring\Lambda^{k}(F)$,
 as follows, for example, by the universal coefficient theorem \cite[Theorem 11.1]{mac1975homology}.
 
 It is easy to see that the linear mapping
 \begin{align*}
  \Theta^m_k : 
  \Kappa^{m}_{k}(F)
  \rightarrow
  \calC_l(\calM_F,\calN_F) \otimes \mathring\Lambda^{k}(F)
 \end{align*}
 which is defined by
 \begin{align*}
  \ext^k_{F,C} \omega
  \mapsto
  C \otimes \omega,
  \quad
  C \in \calT^m \setminus \calU^m,
  \quad
  \omega \in \mathring\Lambda^{k}(F),
 \end{align*}
 is bijective, and that 
 \begin{align*}
  \begin{CD}
   \Kappa^{m}_{k}(F) @>{\Diffver^m_k}>> \Kappa^{m}_{k}(F)
   \\
   @V{\Theta^{m  }_k}VV @V{\Theta^{m-1}_k}VV
   \\
   \calC_{m  }(\calM_F,\calN_F) \otimes \mathring\Lambda^{k}(F)
   @>{ \partial_{m} \otimes \Id }>>
   \calC_{m-1}(\calM_F,\calN_F) \otimes \mathring\Lambda^{k}(F)
  \end{CD}
 \end{align*}
 is a commuting diagram. So $\Theta^m_k$ is an isomorphism of differential complexes from \eqref{seq:sub} to \eqref{seq:dom},
 and thus induces isomorphisms on homology.
 The desired result now follows by the local patch condition.
\end{proof}

\begin{corollary}
 \label{res:trace_complex_is_exact}
 If the geometric decomposition assumption and the local patch condition hold,
 then the Hilbert complex \eqref{eq:trace_complex_extended} is exact.
\end{corollary}


The main result of this section is based on the notion of double complexes in homological algebra.
Let us consider the diagram
\begin{align}
 \label{seq:double_complex}
 \begin{split}
  \begin{CD}
    0                       @>>> \Lambda^{0  }(\calT^{n  }) @>{\Diffhor^{n  }_{0  }}>> \Lambda^{1  }(\calT^{n  }) @>{\Diffhor^{n  }_{1  }}>> \dots
    \\
    @VVV @VVV @VVV @. 
    \\
    \Kappa^{0}(\calT^{n  }) @>>> \Lambda^{0  }_{-1}(\calT^{n  }) @>{ \Diffhor^{n  }_{0  }}>> \Lambda^{1  }_{-1}(\calT^{n  }) @>{ \Diffhor^{n  }_{1  }}>> \dots
    \\
    @V{ -\Diffver^{n  }_{0  }}VV @V{-\Diffver^{n  }_{0  }}VV @V{-\Diffver^{n  }_{1  }}VV @. 
    \\
    \Kappa^{0}(\calT^{n-1}) @>>> \Lambda^{0  }_{-1}(\calT^{n-1}) @>{-\Diffhor^{n-1}_{0  }}>> \Lambda^{1  }_{-1}(\calT^{n-1}) @>{-\Diffhor^{n-1}_{1  }}>> \dots
    \\
    @V{  \Diffver^{n  }_{0  }}VV @V{ \Diffver^{n-1}_{0  }}VV @V{ \Diffver^{n-1}_{1  }}VV @. 
    \\
    \dots @. \dots @. \dots @.
  \end{CD}
 \end{split}
\end{align}
where the left-most horizontal and the top-most vertical arrows denote the respective inclusion mappings.
Note that the choice of signs in \eqref{seq:double_complex} is motivated by \eqref{eq:distributional_differential},
so that the identities 
\begin{align*}
 \Diffhor^{m}_{k+1} \Diffhor^{m}_{k} = 0,
 \quad
 \Diffver^{m-1}_{k} \Diffver^{m}_{k} = 0,
 \quad
 \Diffver^{m}_{k+1} \Diffhor^{m}_{k} - \Diffhor^{m-1}_{k} \Diffver^{m}_{k} = 0
\end{align*}
hold. These identities imply that \eqref{seq:double_complex} constitutes a double complex in the sense of \cite[Chapter 1, §3.5]{gelfand1999homological}.
The previous results on the homology of the horizontal and vertical complexes allow us
to relate the harmonic spaces of two further families of differential complexes:
\begin{align}
 \label{eq:dddf_standard_sequence}
 \begin{CD}
  0
  @>>>
  \Lambda^{0}(\calT^m)
  @>{\Diffhor^m_{0  }}>>
  \dots
  @>{\Diffhor^m_{m-1}}>>
  \Lambda^{m}(\calT^m)
  @>>>
  0,
 \end{CD}
 \\
 \label{eq:dddf_hongkongfu}
 \begin{CD}
  0
  @>>>
  \Lambda^k(\calT^{m})
  @>{\Diffver^{m}_{k}}>>
  \dots
  @>{\Diffver^{k+1}_{k}}>>
  \Lambda^{k}(\calT^{k})
  @>>>
  0.
 \end{CD}
\end{align}
The complex \eqref{eq:dddf_standard_sequence} resembles the standard finite element complex,
and for $m=n$, it coincides with the standard finite element complex of finite element exterior calculus
in our example application.
Similarly, the complex \eqref{eq:dddf_hongkongfu} resembles the chain complex of $\calT$ relative to $\calU$,
and for $k=0$, those two complexes coincide, up to a sign convention, in our example application.

We denote the harmonic spaces of the Hilbert complex \eqref{eq:dddf_standard_sequence} by $\frakH^k(\calT^m)$,
and the harmonic spaces of the Hilbert complex \eqref{eq:dddf_hongkongfu} by $\frakC^k(\calT^m)$.
\begin{align}
 \frakH^k(\calT^m) &:= \left\{ \omega \in \Lambda^{k}(\calT^{m}) \st
 \omega \in \ker \Diffhor^{m}_{k}, \omega \perp \Diffhor^{m}_{k-1}\Lambda^{k}(\calT^{m})
 \right\}
 \\
 \frakC^k(\calT^m) &:= \left\{ \omega \in \Kappa ^{k}(\calT^{m}) \st
 \omega \in \ker \Diffver^{m}_{k}, \omega \perp \Diffver^{m+1}_{k}\Kappa ^{k}(\calT^{m+1})
 \right\}
\end{align}

\begin{theorem}
 Suppose that the rows and columns of the double complex \eqref{seq:double_complex} are exact sequences,
 with the possible exception of the top-most row and the left-most column.
 Then $\frakH^k(\calT^n)$ is isomorphic to $\frakC^0(\calT^{n-k})$ for $0 \leq k \leq n$.
\end{theorem}

\begin{proof}
 This is a standard result in homological algebra;
 see for example Proposition 3.11 of \cite{osborne2000basic},
 Chapter 9.2 of \cite{bott1982differential}
 or Corollary 6.4 of \cite{barr2002acyclic}.
\end{proof}

\begin{corollary}
 \label{res:homology_result_special}
 Suppose that the local exactness condition, the geometric decomposition condition,
 and the local patch condition hold.
 Then we have isomorphisms between harmonic spaces
 \begin{align*}
  \calH_{n-k}(\calT,\calU) \simeq \frakC^0(\calT^{n-k}) \simeq \frakH^k(\calT^n)
 \end{align*}
 for $0 \leq k \leq n$.
\end{corollary}

\begin{example}
 We continue our example application.
 The complex
 \begin{align}
  \begin{CD}
    0
    @>>>
    \Lambda^0(\calT^n)
    @>\cartan^0>>
    \dots
    @>\cartan^{n-1}>>
    \Lambda^n(\calT^n)
    @>>>
    0,
  \end{CD}
 \end{align}
 is a complex of finite element differential forms whose traces on simplices of $\calV$ vanish.
 This means that it is a conforming discretization of the $L^2$ de Rham complex
 with partial tangential boundary conditions along $\Gamma_T$.
 The discrete harmonic forms of the finite element complex are $\frakH^k(\calT)$.
 We have
 \begin{align*}
  \dim \frakH^k(\calT)
  &=
  b_{n-k}(\Omega,\Gamma_N)
  \\&=
  \dim \frakH^{n-k}(\Omega,\Gamma_N,\Gamma_T)
  \\&=
  \dim \frakH^{  k}(\Omega,\Gamma_T,\Gamma_N)
  =
  b_{  k}(\Omega,\Gamma_T)
 \end{align*}
 for $0 \leq k \leq n$.
 This includes the special cases $\Gamma_T = \emptyset$ and $\Gamma_T = \Gamma$, 
 which have been treated earlier in literature \cite{christiansen2008smoothed}.
 \demo
\end{example}

\begin{remark}
 Arnold, Falk and Winther have derived the finite element de Rham cohomology without boundary conditions from the $L^2$ de Rham complex \cite{AFW1};
 Christiansen and Winther have extended this to the case of essential boundary conditions \cite{christiansen2008smoothed}.
 Christiansen has also derived the finite element de Rham cohomology without boundary conditions within the framework of element systems via de Rham mappings \cite{StructPresDisc}.
 With different techniques, we have derived the finite element de Rham cohomology without reference to the $L^2$ de Rham complex.
\end{remark}


\ifx\justbeingincluded\undefined
\end{document}
\fi


\ifx\justbeingincluded\undefined
\input{../global/header_common.tex}

\begin{document}

\fi

\section{Harmonic Forms of Discrete Distributional de Rham Complexes}
\label{sect:distributional_homology}

In this section we introduce discrete distributional de Rham complexes
and construct isomorphisms between their harmonic spaces.
Our main result, Theorem \ref{res:theorema_egregium} below, generalizes Corollary \ref{res:homology_result_special} of the previous section.
In particular, we explicitly construct an isomorphism between $\calH_{n-k}(\calT,\calU)$ and $\frakH^k(\calT^n)$.
%
\\

In this section, we continue to assume that $\calT$ is a finite $n$-dimensional simplicial complex,
and that $\calU$ is a simplicial subcomplex.
\\

Consider again the complex
\begin{align}
 \label{seq:fem_original}
 \begin{CD}
  0
  @>>>
  \Lambda^0(\calT^n)
  @>\cartan^0>>
  \dots
  @>\cartan^{n-1}>>
  \Lambda^n(\calT^n)
  @>>>
  0,
 \end{CD}
\end{align}
which resembles the standard finite element complex.
This complex might be ``redirected'' at any index $k$,
in the sense that we replace $\Lambda^k(\calT^n)$ with $\Lambda^k_{-1}(\calT^n)$,
and continue the complex with the spaces $\Lambda^{k+1}_{-2}(\calT^n)$, $\Lambda^{k+2}_{-3}(\calT^n)$, and so forth,
with the exterior derivative applied in the distributional sense.
\begin{align}
 \label{seq:fem_broken}
 \begin{CD}
  \dots
  @>\cartan^{k-2}>>
  \Lambda_{  }^{k-1}(\calT^n)
  @>\cartan^{k-1}>>
  \Lambda_{-1}^{k  }(\calT^n)
  @>\cartan^{k  }>>
  \Lambda_{-2}^{k+1}(\calT^n)
  @>\cartan^{k+1}>>
  \dots
 \end{CD}
\end{align}
We see that the original complex is already trivially redirected at the $n$-forms,
noting $\Lambda^{n}(\calT^n) = \Lambda^n_{-1}(\calT^n)$.
This is a subcomplex of the complex redirected at the $(n-1)$-forms,
which in turn is a subcomplex of the complex redirected at the $(n-2)$-forms.
We proceed in this manner, until we eventually have a ``maximal'' complex 
that is redirected already at the $0$-forms.
We thus observe a sequence of complexes,
from the original complex over a succession of \emph{discrete distributional de Rham complexes}
to a ``maximal'' complex:
\begin{align}
 \label{seq:fem_maximal}
 \begin{CD}
  0
  @>>>
  \Lambda_{-1}^0(\calT^n)
  @>\cartan^0>>
  \dots
  @>\cartan^{n-1}>>
  \Lambda_{-n-1}^n(\calT^n)
  @>>>
  0.
 \end{CD}
\end{align}
Completely analogous arguments can be applied to the complex
\begin{align}
 \label{seq:chains_original}
 \begin{CD}
  0
  @>>>
  \Kappa^0(\calT^n)
  @>\cartan^{0}>>
  \dots
  @>\cartan^{n-1}>>
  \Kappa^0(\calT^0)
  @>>>
  0.
 \end{CD}
\end{align}
At any index $m$, this complex can be redirected,
\begin{align}
 \label{seq:chains_broken}
 \begin{CD}
  \dots
  @>\cartan^{k-2}>>
  \Kappa_{  }^{0}(\calT^{n-k+1})
  @>\cartan^{k-1}>>
  \Kappa_{-1}^{0}(\calT^{n-k  })
  @>\cartan^{k  }>>
  \Kappa_{-2}^{0}(\calT^{n-k-1})
  @>\cartan^{k+1}>>
  \dots
 \end{CD}
\end{align}
and again we observe a sequence of discrete distributional de Rham complexes.
The maximal example of this second family complexes,
\begin{align}
 \label{seq:chains_maximal}
 \begin{CD}
  0
  @>>>
  \Kappa_{-1  }^0(\calT^n)
  @>\cartan^0>>
  \dots
  @>\cartan^{n-1}>>
  \Kappa_{-n-1}^0(\calT^0)
  @>>>
  0,
 \end{CD}
\end{align}
is, in fact, identical to \eqref{seq:fem_maximal}.
Our goal is to construct isomorphisms between the harmonic spaces of all these complexes.

More generally, we consider discrete distributional de Rham complexes
on the lower dimensional skeletons of $\calT$.
Up to signs, those complexes are described by
\begin{align}
 \label{seq:fem_broken_general}
 \begin{CD}
  \dots
  @>\cartan^{k+n-m-2}>>
  \Lambda_{  }^{k-1}(\calT^m)
  @>\cartan^{k+n-m-1}>>
  \Lambda_{-1}^{k  }(\calT^m)
  @>\cartan^{k+n-m  }>>
  \Lambda_{-2}^{k+1}(\calT^m)
  @>\cartan^{k+n-m+1}>>
  \dots
 \end{CD}
\end{align}
for $0 \leq m \leq n$.
Analogously, albeit not as intuitively, we may consider the discrete distributional de Rham complex
\begin{align}
 \label{seq:chains_broken_general}
 \begin{CD}
  \dots
  @>\cartan^{k+n-m-2}>>
  \Kappa_{  }^{k}(\calT^{m+1})
  @>\cartan^{k+n-m-1}>>
  \Kappa_{-1}^{k}(\calT^{m  })
  @>\cartan^{k+n-m  }>>
  \Kappa_{-2}^{k}(\calT^{m-1})
  @>\cartan^{k+n-m+1}>>
  \dots
 \end{CD}
\end{align}
for $0 \leq k \leq n$.
Whereas the complex \eqref{seq:fem_broken_general}
can be thought of as being derived by skipping the spaces
associated to high-dimensional simplices,
the complex \eqref{seq:chains_broken_general}
can be thought of as skipping spaces with forms of lower degree.

We denote the harmonic spaces of these complexes by
\begin{align}
 \frakH^k_{-b}(\calT^m)
 &:=
 \left\{ \omega \in \Lambda^{k  }_{-b}(\calT^m)
 \st \Difftot^{k+n-m} \omega = 0,
 \;
 \omega \perp \Difftot^{k+n-m-1}\Lambda^{k-1}_{-b+1}(\calT^m)
 \right\},
 \\
 \frakC^k_{-b}(\calT^m)
 &:=
 \left\{ \omega \in \Kappa^{k  }_{-b}(\calT^m)
 \st \Difftot^{k+n-m} \omega = 0,
 \;
 \omega \perp \Difftot^{k+n-m-1}\Kappa^{k}_{-b+1}(\calT^{m+1})
 \right\}
 ,
\end{align}
and use the term \emph{discrete distributional harmonic form} for the elements of these harmonic spaces.
We sometimes write $\frakH^k_{}(\calT^m) = \frakH^k_{0}(\calT^m)$
and $\frakC^k_{}(\calT^m) = \frakC^k_{0}(\calT^m)$ for notational reasons.
We call the elements of these spaces \textit{discrete distributional harmonic forms}.

\begin{remark}
 The results in this section generalize ideas of \cite{BrSchoMax},
 in particular the proofs of their Lemma 3, Theorem 5 and Theorem 7.
 But the distributional complexes in this section can also be 
 related to the double complex of the preceding section.
 We identify the maximal complex \eqref{seq:fem_maximal} / \eqref{seq:chains_maximal}
 as the total complex of the double complex \eqref{seq:double_complex},
 skipping the left-most column and the top-most-row of that diagram.
 The two families of broken complexes, \eqref{seq:fem_broken} and \eqref{seq:chains_broken},
 exemplify the two canonical filtrations of the total complex.
 We refer to \cite{bott1982differential} for more background on notion of homological algebra.
 Although the underlying ideas are similar, our presentation is specifically tailored towards
 finite element analysis and addresses the harmonic spaces of the broken complexes explicitly.
\end{remark}

The above sequences can be defined with the notions of Section \ref{sect:element_systems},
but in order to derive the desired isomorphisms between the harmonic spaces,
we utilize the additional assumptions of Section \ref{sect:homology_theory_theoretical}.
This means that in the sequel, we assume that the local exactness condition, the geometric decomposition condition,
and the local patch condition of $\calT$ relative to $\calU$ hold.
\\


For the construction of the isomorphisms we assume to have right-inverses
of the operators $\Diffhor^{m}_{k}$ and $\Diffver^{m}_{k}$.
This means that we have operators
\begin{align*}
 \Antiver^m_k : \Lambda_{-1}^k(\calT^{m-1}) \rightarrow \Lambda_{-1}^k(\calT^{m}),
 \quad
 \Antihor^m_k : \Lambda_{-1}^{k+1}(\calT^{m}) \rightarrow \Lambda_{-1}^{k}(\calT^{m})
\end{align*}
that satisfy 
\begin{align*}
 \Diffver^m_k = \Diffver^m_k \Antiver^m_k \Diffver^m_k,
 \quad
 \Diffhor^m_k = \Diffhor^m_k \Antihor^m_k \Diffhor^m_k.
\end{align*}
The Moore-Penrose pseudoinverses \cite{noteonpseudoinv} of $\Diffhor^{m}_{k}$ and $\Diffver^{m}_{k}$
are a possible choice for $\Antihor^m_k$ and $\Antiver^m_k$, respectively.
It is notationally convenient to introduce the following operators as well.
We consider, on the one hand,
\begin{align}
 \label{eq:R_operators}
 \begin{array}{rll}
  R_{k,b} &:
  \Lambda_{-b}^k(\calT^n)
  \rightarrow
  \Lambda_{-b}^k(\calT^n),
  &\quad
  \omega
  \mapsto
  \omega - (-1)^{b+1} \Difftot^{k-1} \Antiver^{n-b+2}_{k-b+1} \omega,
  \end{array}
\end{align}
where $2 \leq b \leq k+1$,
and on the other hand,
\begin{align}
 \label{eq:S_operators}
 \begin{array}{rll}
  S_{m,b} &:
  \Kappa^0_{-b}(\calT^m)
  \rightarrow
  \Kappa^0_{-b}(\calT^m),
  &\quad
  \omega
  \mapsto
  \omega + (-1)^{b+n-m} \Difftot^{n-m-1} \Antihor^{m+b-1}_{b-1} \omega,
  \end{array}
\end{align}
where $2 \leq b \leq n-m+1$.

For each of the results below on the broken complexes \eqref{seq:fem_broken} that generalize finite element complexes,
there exists an analogous result on the broken complexes \eqref{seq:chains_broken} that generalize simplicial chain complexes.
In the sequel, we state both results, but only the give the result in former case,
since the proof in latter case is completely analogous.
\\



Our first observation is, loosely speaking, that the images of discrete distributional differential forms
under the discrete distributional exterior derivative
always have preimages that are ``more regular'' than those images. 
The idea follows the following intuition:
Suppose that $\omega \in \Lambda^k_{-1}(\calT^n)$
with $\Difftot^{k} \omega \in \Lambda^{k+1}_{-1}(\calT^{n})$.
Then $\omega \in \Lambda^k(\calT^n)$ by definition,
so $\Difftot^{k} \omega$ even has a preimage in $\Lambda^{k}_{}(\calT^{n})$.
Completely analogously, suppose that $\omega \in \Kappa^0_{-1}(\calT^m)$
with $\Difftot^{n-m} \omega \in \Kappa^{0}_{-1}(\calT^{m-1})$.
Then $\omega \in \Kappa^k(\calT^m)$ by definition.

More generally, a discrete distributional differential forms in $\Lambda^k_{-b}(\calT^n)$
that has a preimage under $\cartan^{k-1}$ in $\Lambda^{k-1}_{-b}(\calT^n)$
already has a preimage in $\Lambda^{k-1}_{-b+1}(\calT^n)$.
But in the general case the construction is more complicated and involves the operators $R_{k,b}$ and $S_{m,b}$,
which, in this sense, can be seen as regularizers of preimages.
%
%
We have

\begin{lemma}
 \label{res:preimage_regular_primal}
 Let $1 \leq b \leq k+1$.
 If $\omega \in \Lambda^k_{-b}(\calT^n)$ with $\Difftot^{k} \omega \in \Lambda^{k+1}_{-b}(\calT^{n})$,
 then $\Difftot^{k} \omega \in \Difftot^{k} \Lambda^{k}_{-b+1}(\calT^{n})$.
 If moreover $2 \leq b$,
 then $\Difftot^{k} R_{k,b} \omega = \Difftot^{k} \omega$ and $R_{k,b} \omega \in \Lambda^k_{-b+1}(\calT^n)$.
\end{lemma}

\begin{proof}
 Let $\omega = \omega^0 + \dots + \omega^{b-1} \in \Lambda^k_{-b}(\calT^n)$
 with $\omega^j \in \Lambda_{-1}^{k-j}(\calT^{n-j})$.
 From $\Difftot^{k} \omega \in \Lambda^{k+1}_{-b}(\calT^{n})$
 we see $\Diffver^{n-b+1}_{k-b+1} \omega^{b-1} = 0$.
 So $\Diffver^{n-b+2}_{k-b+1} \Antiver^{n-b+2}_{k-b+1} \omega^{b-1} = \omega^{b-1}$
 holds by Corollary \ref{res:trace_complex_is_exact}.
 We conclude that
 \begin{align*}
  R_{k,b} \omega
  &=
  \omega^0 + \dots + \omega^{b-1}
  - (-1)^{b-1} \Difftot^{k-1} \Antiver^{n-b+2}_{k-b+1} \omega^{b-1}
  \\&=
  \omega^0 + \dots + \omega^{b-1}
  + \Diffhor^{n-b+2}_{k-b+1} \Antiver^{n-b+2}_{k-b+1} \omega^{b-1}
  - \Diffver^{n-b+2}_{k-b+1} \Antiver^{n-b+2}_{k-b+1} \omega^{b-1}
  \\&=
  \omega^0 + \dots + \omega^{b-2}
  + \Diffhor^{n-b+2}_{k-b+1} \Antiver^{n-b+2}_{k-b+1} \omega^{b-1},
 \end{align*}
 so $R_{k,b} \omega \in \Lambda^k_{-b+1}(\calT^n)$. Furthermore,
 \begin{align*}
  \Difftot^{k} R_{k,b} \omega
  =
  \Difftot^{k} \omega - (-1)^{b-1} \Difftot^{k} \Difftot^{k-1} \Antiver^{n-b+2}_{k-b+1} \omega^{b-1}
  =
  \Difftot^{k} \omega.
 \end{align*}
 This completes the proof.
\end{proof}

\begin{lemma}
 \label{res:preimage_regular_analog}
 Let $1 \leq b \leq n-m+1$.
 If $\omega \in \Kappa^0_{-b}(\calT^m)$ with $\Difftot^{n-m} \omega \in \Kappa^{0}_{-b}(\calT^{m-1})$,
 then $\Difftot^{n-m} \omega \in \Difftot^{n-m} \Kappa^{0}_{-b+1}(\calT^{m-1})$.
 If moreover $2 \leq b$, 
 then $\Difftot^{n-m} S_{m,b} \omega = \Difftot^{n-m} \omega$ and $S_{m,b} \omega \in \Kappa^0_{-b+1}(\calT^m)$.
\end{lemma}

Another auxiliary result restricts the class of discrete distributional differential forms
that are candidates for being discrete distributional harmonic forms.
The result implies that an element of $\frakH^{k}_{-b+1}(\calT^{n})$ is not contained in $\frakH^{k}_{-b}(\calT^{n})$.
Analogously, an element of $\frakC^{0}_{-b+1}(\calT^{m})$ is not contained in $\frakC^{0}_{-b}(\calT^{m})$.

\begin{lemma}
 \label{res:harmonic_forms_irregular_primal}
 Let $2 \leq b \leq k+1$.
 If $\omega \in \Lambda_{-b+1}^{k}(\calT^{n})$
 with $\cartan^{k}\omega = 0$,
 then $\omega$ is not orthogonal to $\cartan^{k-1} \Lambda_{-b+1}^{k-1}(\calT^{n})$.
\end{lemma}

\begin{proof}
 Suppose that $\omega = \omega^{0} + \dots + \omega^{b-2} \in \Lambda^{k}_{-b}(\calT^{n})$
 with $\omega^{j} \in \Lambda^{k-j}(\calT^{n-j})$ and $\cartan^{k} \omega = 0$.
 We then set $\xi^{0} := \Antihor^{n}_{k-1} \omega^{0}$,
 and define recursively $\xi^j \in \Lambda^{k-j-1}(\calT^{n-j})$ by
 \begin{align*}
  \xi^{j}
  &:=
  (-1)^{j} \Antihor^{n-j}_{k-j-1} \left(
  \omega^{j} - (-1)^{j} \Diffver^{n-j+1}_{k-j} \xi^{j-1} \right),
  \quad 1 \leq j \leq b-2.
 \end{align*}
 We clearly have $\cartan^{k-1} \xi^{0} = \omega^{0} - \Diffver^{n}_{k-1} \xi^{0}$, 
 since $\Diffhor^{m}_{k} \omega^{0} = 0$ and the local exactness condition holds.
 
 Now assume that we have already shown
 \begin{align*}
  \cartan^{k-1} ( \xi^{0} + \dots + \xi^{j} )
  =
  \omega^{0} + \dots + \omega^{j} + (-1)^{j+1} \Diffver^{n-j}_{k-j-1} \xi^{j}
 \end{align*}
 for $j < b-2$. This equation implies in particular that
 \begin{align*}
  (-1)^{j} \Diffhor^{n-j}_{k-j-1} \xi^{j}
  =
  \omega^{j} - (-1)^{j} \Diffver^{n-j+1}_{k-j} \xi^{j-1},
 \end{align*}
 where we set $\xi^{-1} = 0$ for notational reasons.
 So we find that
 \begin{align*}
  \Diffhor^{n-j-1}_{k-j-1} \Diffver^{n-j}_{k-j-1} \xi^{j}
  &=
  \Diffver^{n-j}_{k-j} \Diffhor^{n-j}_{k-j-1} \xi^{j}
  \\&=
  (-1)^{j} \Diffver^{n-j}_{k-j} \omega^{j}
  - \Diffver^{n-j}_{k-j} \Diffver^{n-j+1}_{k-j} \xi^{j-1}
  \\&=
  (-1)^{j} \Diffver^{n-j}_{k-j} \omega^{j},
 \end{align*}
 and calculate, using $\cartan^{k}\omega = 0$, that
 \begin{align*}
  &\quad
  \Diffhor^{n-j-1}_{k-j-1} \left( \omega^{j+1} - (-1)^{j+1} \Diffver^{n-j}_{k-j-1} \xi^{j} \right)
  \\&=
  \Diffhor^{n-j-1}_{k-j-1} \omega^{j+1}
  -
  (-1)^{j+1} \Diffhor^{n-j-1}_{k-j-1} \Diffver^{n-j}_{k-j-1} \xi^{j}
  \\&=
  \Diffhor^{n-j-1}_{k-j-1} \omega^{j+1}
  + 
  \Diffver^{n-j}_{k-j} \omega^{j}
  \\&=
  0.
 \end{align*}
 The local exactness condition implies that
 \begin{align*}
  \Diffhor^{n-j-1}_{k-j-2} \xi^{j+1}
  =
  (-1)^{j+1} \omega^{j+1} - \Diffver^{n-j}_{k-j-1} \xi^{j}
  . 
 \end{align*}
 We thus find
 \begin{align*}
  &\quad
  \cartan^{k-1} ( \xi^{0} + \dots + \xi^{j} + \xi^{j+1} )
  \\&=
  \omega^{0} + \dots + \omega^{j} 
  + (-1)^{j+1} \Diffver^{n-j}_{k-j-1} \xi^{j}
  \\&\quad\quad
  + (-1)^{j+1} \Diffhor^{n-j-1}_{k-j-2} \xi^{j+1}
  + (-1)^{j  } \Diffver^{n-j-1}_{k-j-2} \xi^{j+1}
  \\&=
  \omega^{0} + \dots + \omega^{j} + (-1)^{j+1} \Diffver^{n-j}_{k-j-1} \xi^{j}
  \\&\quad\quad
  + \omega^{j+1} - (-1)^{j+1} \Diffver^{n-j}_{k-j-1} \xi^{j}
  + (-1)^{j  } \Diffver^{n-j-1}_{k-j-2} \xi^{j+1}
  \\&=
  \omega^{0} + \dots + \omega^{j}
  + \omega^{j+1}
  + (-1)^{j  } \Diffver^{n-j-1}_{k-j-2} \xi^{j+1}
  .
 \end{align*}
 So eventually we see that
 \begin{align*}
  \cartan^{k-1} ( \xi^{0} + \dots + \xi^{b-2} )
  =
  \omega^{0} + \dots + \omega^{b-2} + (-1)^{b-1} \Diffver^{n-b+2}_{k-b+1} \xi^{b-2},
 \end{align*}
 from which we deduce
 \begin{align*}
  \left\langle \cartan^{k-1} ( \xi^{0} + \dots + \xi^{b-2} ), \omega \right\rangle = \| \omega \|^{2},
 \end{align*}
 proving the claim.
\end{proof}

\begin{lemma}
 \label{res:harmonic_forms_irregular_analog}
 Let $2 \leq b \leq n-m+1$.
 If $\omega \in \Kappa_{-b+1}^{0}(\calT^{m})$
 with $\cartan^{n-m}\omega = 0$,
 then $\omega$ is not orthogonal to $\cartan^{n-m-1} \Kappa_{-b+1}^{0}(\calT^{m+1})$.
\end{lemma}


Our goal is to construct the discrete distributional harmonic forms of the distributional complexes,
and to find isomorphisms between the harmonic spaces.
To begin with, the harmonic forms of the spaces $\Lambda^k_{-1}(\calT^n)$ and $\Kappa^k_{-1}(\calT^n)$ are easily described:

\begin{lemma}
 \label{res:initial_harmonics_primal}
 We have $\frakH^k(\calT^n) = \frakH^k_{-1}(\calT^n)$
 for $0 \leq k \leq n$.
\end{lemma}

\begin{proof}
 We know that $\omega \in \Lambda_{-1}^k(\calT^n)$ with $\Difftot^k \omega = 0$,
 is equivalent to $\omega \in \Lambda^k(\calT^n)$ with $\Difftot^k \omega = 0$.
 The equality $\frakH^k(\calT^n) = \frakH^k_{-1}(\calT^n)$ now follows from definitions.
\end{proof}

\begin{lemma}
 \label{res:initial_harmonics_analog}
 We have $\frakC^0(\calT^m) = \frakC^0_{-1}(\calT^m)$
 for $0 \leq m \leq n$.
\end{lemma}

The harmonic spaces $\frakH^{k}_{-b}(\calT^{n})$ and $\frakC^{0}_{-b}(\calT^{m})$ for $b \geq 2$
are constructed in a recursive manner.

\begin{lemma}
 \label{res:recursive_harmonics_primal}
 Suppose that $2 \leq b \leq k+1$,
 and let $Q^{k}_{b}$ be the orthogonal projection onto the kernel
 of the operator $\cartan^{k} : \Lambda^{k}_{-b}(\calT^{n}) \rightarrow \Lambda^{k+1}_{-b+1}(\calT^{n})$.
 Then the operator $Q^{k}_{b} R_{k,b}^{\ast}$ acts as an isomorphism
 from $\frakH^{k}_{-b+1}(\calT^{n})$ to $\frakH^{k}_{-b}(\calT^{n})$.
\end{lemma}

\begin{proof}
 Let $\omega = \omega^{0} + \dots + \omega^{b-1} \in \Lambda^{k}_{-b}(\calT^{n})$
 with $\omega^{j} \in \Lambda_{-1}^{k-j}(\calT^{n-j})$.
 We have by construction that $\omega - R_{k,b} \omega \in \cartan^{k-1}\Lambda^{k-1}_{-b+1}(\calT^{n})$.
 This implies in particular that
 \begin{align*}
  \cartan^{k} \omega = 0
  &\quad\Longleftrightarrow\quad
  \cartan^{k} R_{k,b} \omega = 0,
  \\
  \omega \in \cartan^{k-1}\Lambda^{k-1}_{-b+1}(\calT^{n})
  &\quad\Longleftrightarrow\quad
  R_{k,b} \omega \in \cartan^{k-1}\Lambda^{k-1}_{-b+1}(\calT^{n}).
 \end{align*}
 From the last equivalence and the abstract Hodge decomposition, we conclude that
 \begin{align*}
  \cartan^{k} \omega = 0, \quad \omega \perp \frakH^{k}_{-b}(\calT^{n})
  &\quad\Longleftrightarrow\quad
  \cartan^{k} R_{k,b} \omega = 0, \quad R_{k,b} \omega \perp \frakH^{k}_{-b}(\calT^{n})
  .
 \end{align*}
 Now $\cartan^{k} R_{k,b} \omega = 0$ implies that $R_{k,b} \omega \in \Lambda^{k}_{-b+1}(\calT^{n})$
 by Lemma \ref{res:preimage_regular_primal}. So
 \begin{align*}
  R_{k,b} \omega \in \cartan^{k-1}\Lambda^{k-1}_{-b+1}(\calT^{n})
  \quad\Longleftrightarrow\quad
  R_{k,b} \omega \in \cartan^{k-1}\Lambda^{k-1}_{-b+2}(\calT^{n})
 \end{align*}
 by Lemma \ref{res:preimage_regular_primal} again.
 We derive from this that
 \begin{align*}
  &\quad
  R_{k,b} \omega \in \cartan^{k-1}\Lambda^{k-1}_{-b+2}(\calT^{n})
  \\ \Longleftrightarrow&\quad
  \cartan^{k} R_{k,b} \omega = 0, \quad R_{k,b} \omega \perp \frakH^{k-1}_{-b+1}(\calT^{n})
  \\ \Longleftrightarrow&\quad
  \cartan^{k} \omega = 0, \quad \omega \perp R_{k,b}^{\ast} \frakH^{k-1}_{-b+1}(\calT^{n}).
 \end{align*}
 We conclude that the projection of $R_{k,b}^{\ast} \frakH^{k}_{-b+1}(\calT^{n})$
 onto $\ker \cartan^{k} \cap \Lambda^{k}_{-b}(\calT^n)$
 equals $\frakH^{k}_{-b}(\calT^{n})$.
 Furthermore, we observe for $p \in \frakH^{k}_{-b+1}(\calT^{n})$ that
 \begin{align*}
    \langle p, Q^{k}_{b} R_{k,b}^{\ast} p \rangle
  = \langle p, R_{k,b}^{\ast} p \rangle
  = \langle R_{k,b} p, p \rangle
  = \langle p, p \rangle.
 \end{align*}
 This is a consequence of Lemma \ref{res:harmonic_forms_irregular_primal}.
 We conclude that $Q^{k}_{b} R_{k,b}^{\ast}$ defines an isomorphism
 from $\frakH^{k}_{-b+1}(\calT^{n})$ onto $\frakH^{k}_{-b}(\calT^{n})$.
\end{proof}

\begin{lemma}
 \label{res:recursive_harmonics_analog}
 Suppose that $2 \leq b \leq n-m+1$,
 and let $Q^m_b$ be the orthogonal projection onto the kernel
 of the operator $\cartan^{n-m} : \Kappa^{0}_{-b}(\calT^{m}) \rightarrow \Kappa^{0}_{-b+1}(\calT^{m-1})$.
 Then the operator $Q^m_b S_{m,b}^{\ast}$ acts as an isomorphism
 from $\frakC^{0}_{-b+1}(\calT^{m})$ to $\frakC^{0}_{-b}(\calT^{m})$.
\end{lemma}

\begin{remark}
 If $p \in \frakH_{-b}^{k}(\calT^{n})$, then generally $R_{k,b+1}^{\ast} p \notin \frakH_{-b}^{k}(\calT^{n})$.
 However, in the special case $k=n$ the orthogonal projection
 is redundant because $\cartan^{n} \Lambda_{-b}^{n}(\calT^{n}) = 0$.
 
 The requirement of an orthogonal projection in the construction of the discrete distributional harmonic forms
 seems conceptually unsatisfying.
 However, one can see that, if we leave out this projection operators at every stage,
 the resulting construction produces at each stage a space of discrete distributional harmonic forms
 whose projection onto the discrete distributional harmonic forms is surjective.
 
 Analogous observations apply to the operators $S_{m,b}$.
\end{remark}


The main result of this contribution is now evident from Lemmas
\ref{res:initial_harmonics_primal} and \ref{res:initial_harmonics_analog}, 
and the repeated application of Lemmas \ref{res:recursive_harmonics_primal} and \ref{res:recursive_harmonics_analog}.
It generalizes of Corollary \ref{res:homology_result_special}.

\begin{theorem}
 \label{res:theorema_egregium}
 Under the assumptions of this section, 
 we have isomorphisms between harmonic spaces:
 \begin{align*}
  \calH^m(\calT,\calU)
  \simeq
  \frakC^0(\calT^{n-k})
  &=
  \frakC^0_{  -1}(\calT^{n-k})
  \simeq
  \dots
  \simeq
  \frakC^0_{-k-1}(\calT^{n-k})
  \\&=
  \frakH^k_{-k-1}(\calT^{n})
  \simeq
  \dots
  \simeq
  \frakH^k_{  -1}(\calT^n)
  =
  \frakH^k(\calT^n)
 \end{align*}
 for $0 \leq k \leq n$.
 \qed
\end{theorem}



We have studied the harmonic spaces $\frakH^k_{-b}(\calT^n)$ and $\frakC^0_{-b}(\calT^{n-k})$.
The harmonic spaces $\frakH^{k}_{-b}(\calT^{n-1})$ can be studied in a similar manner.
One merely replaces $\calT^{n}$ by $\calT^{[n-1]} \setminus \calU^{n-1}$ and $\calU^{}$ by $\calU^{[n-2]}$,
and uses the arguments of the previous and this section.
Note that the local patch condition is then satisfied; see Remark \ref{res:local_patch_condition_skeleton} above.
Repeating this idea for $0 \leq m < n$, we obtain isomorphisms
 \begin{align*}
  \frakC^0(\calT^{m	-k})
  &=
  \frakC^0_{  -1}(\calT^{m-k})
  \simeq
  \dots
  \simeq
  \frakC^0_{-k-1}(\calT^{m-k})
  \\&=
  \frakH^k_{-k-1}(\calT^{m})
  \simeq
  \dots
  \simeq
  \frakH^k_{  -1}(\calT^m)
  =
  \frakH^k(\calT^m),
  \quad 0 < k \leq m.
 \end{align*}
However, the case $k=0$ must be treated differently,
since in general 
\begin{align*}
 \calH_{n-1}(\calT^{[n-1]}\setminus\calU^{n-1},\calU^{n-2}) \neq \calH_{n-1}(\calT,\calU).
\end{align*}
We can use that the harmonic space in the discrete distributional $0$-forms
for the discrete distributional de Rham complex over the $m$-skeleton of $\calT$ is the direct orthogonal sum
of $\frakC^{0}(\calT^{m})$ and $\Diffver^{m+1}_{0} \Kappa^{0}(\calT^{m+1})$.
This immediately obvious, since that complex can be derived by ``truncation''
of a discrete distributional de Rham complex of the form \eqref{seq:chains_broken_general}.
This allows us to relate the harmonic spaces of discrete distributional de Rham complexes
defined over different skeletons of $\calT$.
The following result is of ancillary interest.

\begin{lemma}
 The orthogonal projection 
 from $\Lambda^{k}_{-2}(\calT^{n})$ onto the kernel
 of $\cartan^{k} : \Lambda^{k-1}(\calT^{n-1}) \rightarrow \Lambda^{k}(\calT^{n-1})$
 maps $\frakH^{k}_{-2}(\calT^{n})$ isomorphically onto $\frakH^{k-1}(\calT^{n-1})$
 for $k \geq 2$.
\end{lemma}

\begin{proof}
 Suppose that $\omega \in \cartan^{k+1}\Lambda^{k}(\calT^{n-1})$,
 with $\xi \in \Lambda^{k}(\calT^{n-1})$ such that $\cartan^{k+1} \xi = \omega$.
 Then $R_{k+1,2} \xi \in \Lambda^{k+1}_{-1}(\calT^{n})$ with $\cartan^{k+1} R_{k+1,2} \xi = \omega$
 as follows from Lemma \ref{res:preimage_regular_primal}.
 We conclude that
 \begin{align*}
  &
  \cartan^{k+1} \omega = 0, \quad \omega \perp \frakH^{k}(\calT^{n-1})
  \\ \Longleftrightarrow\quad&
  \omega \in \cartan^{k+1}\Lambda^{k}(\calT^{n-1})
  \\ \Longleftrightarrow\quad&
  \omega \in \cartan^{k+1}\Lambda_{-1}^{k+1}(\calT^{n})
  \\ \Longleftrightarrow\quad&
  \cartan^{k+1} \omega = 0, \quad \omega \perp \frakH_{-2}^{k+1}(\calT^{n})
 \end{align*}
 Thus we see that the orthogonal projection
 of $\frakH_{-2}^{k+1}(\calT^{n})$
 onto $\Lambda^{k+1}(\calT^{n-1})$
 yields $\frakH^{k}(\calT^{n-1})$.
 Furthermore, from Lemma \ref{res:harmonic_forms_irregular_primal} we conclude
 that this mapping is not only onto, but also one-to-one.
\end{proof}

The harmonic spaces of the complexes \eqref{seq:chains_broken_general}
can be analyzed in an analogous manner. 
We do not describe this in detail, since the differences are mostly notational.

\ifx\justbeingincluded\undefined
\end{document}
\fi



\ifx\justbeingincluded\undefined
\input{../global/header_common.tex}

\begin{document}

\fi

\section{Conclusions}

Complexes of discrete distributional differential forms have been introduced into finite element exterior calculus.
We have analyzed their homology theory in this contribution.
We will analyze Poincar\'e-Friedrichs inequalities in a subsequent contribution.
\\

Applications in a posteriori error estimation motivate this research,
but distributional finite elements appear in other facets of computational partial differential equations as well.
Most prominently, these include non-conforming methods like discontinuous Galerkin finite element methods and finite volume methods.
Furthermore, a distributional elasticity complex in three dimensions appears in the context of Regge calculus \cite{christiansen2011linearization}.
%
Our example application concerns the $L^2$ de Rham complex on a triangulated manifold,
but the discrete theory applies to a larger class of triangulations.
For example, such instances of non-manifold triangulations appear in the numerical treatment of the Electric Field Integral Equation;
see also \cite[Section 5.2]{christiansen2004characterization} for more details.
\\

The idea of distributional de Rham complexes emerged several times in analysis,
and in finite element analysis within at least one other context:
as observed in \cite{StructPresDisc} and \cite{christiansen2013high} within the framework of element systems,
the degrees of freedom in finite element exterior calculus constitute a differential complex by themselves.
The complexes of discrete distributional differential forms emerge in that context again;
for example, the complex
\begin{align*}
 \begin{CD}
  \calP^{-}_{1}\Lambda^{0  }_{-1  }(\calT^{n})
  @>{\cartan^{  0}}>>
  \calP^{-}_{1}\Lambda^{  1}_{-2  }(\calT^{n})
  @>{\cartan^{  1}}>>
  \dots
  @>{\cartan^{n-1}}>>
  \calP^{-}_{1}\Lambda^{  n}_{-n-1}(\calT^{n})
 \end{CD}
\end{align*}
is isomorphic to the complex of degrees of freedom
of the finite element complex
\begin{align*}
 \begin{CD}
  \calP^{ }_{1  }\Lambda^{n  }(\calT)
  @<{\cartan^{n-1}}<<
  \calP^{ }_{2  }\Lambda^{n-1}(\calT)
  @<{\cartan^{n-2}}<<
  \dots
  @<{\cartan^{0}}<<
  \calP^{ }_{n+1}\Lambda^{  0}(\calT).
 \end{CD}
\end{align*}
Further exploration of this relation will contribute to finite element theory in general.
\\

We have considered only finite-dimensional spaces of distributional differential forms.
The ideas of this contribution can be possibly be extended to Hilbert complexes of in\-fi\-ni\-te-di\-men\-sio\-nal spaces.
Such a generalization might be contributive to the convergence analysis of finite element methods.

\section*{Acknowledgments}

The author would like to thank 
S\"oren Bartels for drawing the author's attention to \cite{BrSchoMax} and supervising the diploma thesis from which this work evolved,
and Snorre Christiansen for productive discussions and sharing his notes on double complexes with the author.
Helpful remarks by Jeonghun Lee are appreciated.

\ifx\justbeingincluded\undefined
\end{document}
\fi


\end{document}